\documentclass[a4paper]{article}

\usepackage[margin=2.5cm]{geometry} 

\usepackage{amsmath}
\usepackage{amsthm}
\usepackage{amssymb}

\usepackage[utf8]{inputenc}
\usepackage{hyperref}
\hypersetup{
	unicode,
	colorlinks,
	urlcolor=blue, 
	linkcolor=blue, 
	citecolor=blue,
	pdfproducer={LaTeX},
	pdfcreator={pdflatex}
}

\usepackage{listings}

\usepackage{pythonhighlight}

\usepackage{enumitem}
\setlist{nolistsep}
\usepackage{xcolor}
\definecolor{lightergray}{rgb}{0.9,0.9,0.9}

\usepackage{xspace}

\usepackage{algorithm}
\usepackage{algpseudocode}
\usepackage[capitalise]{cleveref}
\usepackage[most]{tcolorbox}
\tcbuselibrary{listings,breakable}

\usepackage[frozencache=true,cachedir=.]{minted}
\definecolor{darkspringgreen}{rgb}{0.09, 0.45, 0.27}

\usepackage{tikz}
\usetikzlibrary{arrows,shapes,positioning,shadows,trees}

\tikzset{
  basic/.style  = {draw, 
                   text width=6cm, 
                   font=\sffamily, 
                   rectangle},
  root/.style   = {basic, 
                   rounded corners=3pt, 
                   thin, 
                   align=center, 
                   fill=blue!50,
                   minimum height=0.7cm},
  level 2/.style = {basic, 
                    rounded corners=3pt, 
                    thin,align=center, 
                    fill=blue!23,
                    text width=1.85cm,
                    minimum height=0.5cm},
  level 3/.style = {basic, 
                    thin, 
                    align=left, 
                    fill=blue!5, 
                    text width=1.75cm}
}

\newcommand{\im}{u}

\newcommand{\sysmat}{A}
\newcommand{\sino}{b}

\usepackage{soul}
\sethlcolor{lightergray}

\newcommand{\code}[1]{\pyth{#1}\xspace}

\DeclareMathOperator*{\argmin}{arg\,min}

\newcommand{\fun}{\code{Function}}
\newcommand{\op}{\code{Operator}}
\newcommand{\alg}{\code{Algorithm}}
\newcommand{\proc}{\code{Processor}}
\newcommand{\cilversion}{21.0}
\newcommand{\module}{\textbf}

\usepackage{scalerel}
\usetikzlibrary{svg.path}
 
\definecolor{orcidlogocol}{HTML}{A6CE39}
\tikzset{
  orcidlogo/.pic={
    \fill[orcidlogocol] svg{M256,128c0,70.7-57.3,128-128,128C57.3,256,0,198.7,0,128C0,57.3,57.3,0,128,0C198.7,0,256,57.3,256,128z};
    \fill[white] svg{M86.3,186.2H70.9V79.1h15.4v48.4V186.2z}
                 svg{M108.9,79.1h41.6c39.6,0,57,28.3,57,53.6c0,27.5-21.5,53.6-56.8,53.6h-41.8V79.1z M124.3,172.4h24.5c34.9,0,42.9-26.5,42.9-39.7c0-21.5-13.7-39.7-43.7-39.7h-23.7V172.4z}
                 svg{M88.7,56.8c0,5.5-4.5,10.1-10.1,10.1c-5.6,0-10.1-4.6-10.1-10.1c0-5.6,4.5-10.1,10.1-10.1C84.2,46.7,88.7,51.3,88.7,56.8z};
  }
}

\newcommand\orcidicon[1]{\href{https://orcid.org/#1}{\mbox{\scalerel*{
\begin{tikzpicture}[yscale=-1,transform shape]
\pic{orcidlogo};
\end{tikzpicture}
}{|}}}}

\title{Core Imaging Library -- Part I: a versatile Python framework for tomographic imaging}
\author{
Jakob S. J\o{}rgensen\,\orcidicon{0000-0001-9114-754X}$^{1,8,*}$ \and 
Evelina Ametova\,\orcidicon{0000-0002-8867-3001}$^{2,5}$ \and 
Genoveva Burca\,\orcidicon{0000-0001-6867-9628}$^{3,8}$ \and 
Gemma Fardell\,\orcidicon{0000-0003-2388-5211}$^{4}$ \and
Evangelos Papoutsellis\,\orcidicon{0000-0002-1820-9916}$^{4,5}$ \and
Edoardo Pasca\,\orcidicon{0000-0001-6957-2160}$^{4}$ \and
Kris Thielemans\,\orcidicon{0000-0002-5514-199X}$^{6}$ \and
Martin Turner\,\orcidicon{0000-0003-0117-8049}$^{7}$ \and
Ryan Warr\,\orcidicon{0000-0002-7904-0560}$^{5}$ \and
William R. B. Lionheart\,\orcidicon{0000-0003-0971-4678}$^{8}$ \and
Philip J. Withers\,\orcidicon{0000-0002-1946-5647}$^{5}$}

\date{
	$^{1}$Department of Applied Mathematics and Computer Science, Technical University of Denmark, Richard Petersens Plads, Building 324, 2800 Kgs. Lyngby, Denmark\\
	$^{2}$Laboratory for Applications of Synchrotron Radiation, Karlsruhe Institute of Technology, Kaiserstr. 12, 76131, Karlsruhe, Germany \\
	$^{3}$ISIS Neutron and Muon Source, STFC, UKRI, Rutherford Appleton Laboratory, Harwell Campus, Didcot OX11 0QX, United Kingdom \\
    $^{4}$Scientific Computing Department, STFC,  UKRI,  Rutherford Appleton Laboratory, Harwell Campus, Didcot OX11 0QX, United Kingdom\\ 
    $^{5}$Henry Royce Institute, Department of Materials, The University of Manchester, Oxford Road, Manchester, M13 9PL, United Kingdom\\
    $^{6}$Institute of Nuclear Medicine, University College London, London, United Kingdom\\
    $^{7}$Research IT Services, The University of Manchester, Oxford Road, Manchester M13 9PL, United Kingdom\\
$^{8}$Department of Mathematics, The University of Manchester, Oxford Road, Manchester, M13 9PL, United Kingdom\\
    $^*$Corresponding author: \url{jakj@dtu.dk}}

\begin{document}
	\maketitle
	
	\begin{abstract}
We present the Core Imaging Library (CIL), an open-source Python framework for tomographic imaging with particular emphasis on reconstruction of challenging datasets. Conventional filtered back-projection reconstruction tends to be insufficient for highly noisy, incomplete, non-standard or multi-channel data arising for example in dynamic, spectral and in situ tomography. CIL provides an extensive modular optimisation framework for prototyping reconstruction methods including sparsity and total variation regularisation, as well as tools for loading, preprocessing and visualising tomographic data.
The capabilities of CIL are demonstrated on a synchrotron example dataset and three challenging cases spanning golden-ratio neutron tomography, cone-beam X-ray laminography and positron emission tomography.
	\end{abstract}

 \section{Introduction} \label{sec:introduction}
 
It is an exciting time for computed tomography (CT): existing imaging techniques are being pushed beyond current limits on resolution, speed and dose, while new ones are being continually developed \cite{roadmap2008}. Driving forces include higher-intensity X-ray sources and photon-counting detectors enabling respectively fast time-resolved and energy-resolved imaging. In situ imaging of evolving processes and unconventional sample geometries such as laterally extended samples are also areas of great interest. Similar trends are seen across other imaging areas, including transmission electron microscopy (TEM), positron emission tomography (PET), magnetic resonance imaging (MRI), and neutron imaging, as well as joint or multi-contrast imaging combining several such modalities.

Critical in CT imaging is the reconstruction step where the raw measured data is computationally combined into reconstructed volume (or higher-dimensional) data sets. 
Existing reconstruction software such as proprietary programs on commercial scanners are often optimised for conventional, high quality data sets, relying on filtered back projection (FBP) type reconstruction methods \cite{Pan_2009}.  Noisy, incomplete, non-standard or multi-channel data will generally be poorly supported or not at all.

In recent years, numerous reconstruction methods for new imaging techniques have been developed. In particular, iterative reconstruction methods based on solving suitable optimisation problems, such as sparsity and total variation regularisation, have been applied with great success to improve reconstruction quality in challenging cases \cite{iterativereview}. This however is highly specialised and time-consuming work that is rarely deployed for routine use. The result  is a lack of suitable reconstruction software, severely limiting the full exploitation of new imaging opportunities. 

This article presents the Core Imaging Library (CIL) -- a versatile open-source Python library for processing and reconstruction of challenging tomographic imaging data. CIL is developed by the Collaborative Computational Project in Tomographic Imaging (CCPi) network and is available from \url{https://www.ccpi.ac.uk/CIL} as well as from \cite{CILReleases}, with documentation, installation instructions and numerous demos.

Many software libraries for tomographic image processing already exist, such as TomoPy \cite{Gursoy:tomopy}, ASTRA \cite{vanAarle:16}, TIGRE \cite{Biguri_2016}, Savu \cite{Atwood_2015}, AIR Tools II \cite{airtoolsii2018}, and CASToR \cite{Merlin_2018}. Similarly, many MATLAB and Python toolboxes exist for specifying and solving optimisation problems relevant in imaging, including FOM \cite{BeckFOM2019}, GlobalBioIm \cite{Soubies_2019}, ODL \cite{jonas_adler_2018_1442734}, ProxImaL \cite{Proximal}, and TFOCS \cite{TFOCS}. 

CIL aims to combine the best of the two worlds of tomography and optimisation software in a single easy-to-use, highly modular and configurable Python library. 
Particular emphasis is on enabling a variety of regularised reconstruction methods within a ``plug and play'' structure in which different data fidelities, regularisers, constraints and algorithms can be easily selected and combined. The intention is that users will be able to use the existing reconstruction methods provided, or prototype their own, to deal with noisy, incomplete, non-standard and multi-channel tomographic data sets for which conventional FBP type methods and proprietary software fail to produce satisfactory results.
In addition to reconstruction, CIL supplies tools for loading, preprocessing, visualising and exporting data for subsequent analysis and visual exploration.

CIL easily connects with other libraries to further combine and expand capabilities; we describe CIL plugins for ASTRA \cite{vanAarle:16}, TIGRE \cite{Biguri_2016} and the CCPi-Regularisation (CCPi-RGL) toolkit \cite{Kazantsev2019RGLTK}, as well as interoperability with the Synergistic Image Reconstruction Framework (SIRF) \cite{SIRF2020} enabling PET and MR reconstruction using CIL.

We envision that in particular two types of researchers might find CIL useful:
\begin{itemize}
\item Applied mathematicians and computational scientists can use existing mathematical building blocks and the modular design of CIL to rapidly implement and experiment with new reconstruction algorithms and compare them against existing state-of-the-art methods. They can easily run controlled simulation studies with test phantoms and within the same framework transition into demonstrations on real CT data. 
\item CT experimentalists will be able to load and pre-process their standard or non-standard data sets and reconstruct them using a range of different state-of-the-art reconstruction algorithms. In this way they can experiment with, and assess the efficacy of, different methods for compensating for poor data quality or handle novel imaging modalities in relation to whatever specific imaging task they are interested in.
\end{itemize}
CIL includes a number of standard test images as well as demonstration data and scripts that make it easy for users of both groups to get started using CIL for tomographic imaging. These are described in the CIL documentation and we also highlight that all data and code for the experiments presented here are available as described under Data Accessibility.

This paper describes the core functionality of CIL and demonstrates its capabilities using an illustrative running example, followed by three specialised exemplar case studies. 
\Cref{sec:software} gives an overview of CIL and describes the functionality of all the main modules. \Cref{sec:optimisation} focuses on the optimisation module used to specify and solve reconstruction problems. \Cref{sec:examples} presents the three exemplar cases, before a discussion and outlook are provided in \Cref{sec:discussion}.
Multi-channel functionality (e.g. for dynamic and spectral CT) is presented in the part II paper \cite{CIL2} and a use case of CIL for PET/MR motion compensation is given in \cite{SIRFspecialissue}; further applications of CIL in hyperspectral X-ray and neutron tomography are presented in \cite{ryan} and \cite{Evelina}.

 \section{Overview of CIL} \label{sec:software}
\begin{figure}[t]
\centering
\includegraphics[clip,trim=0 1.0cm 0 0, width=\linewidth]{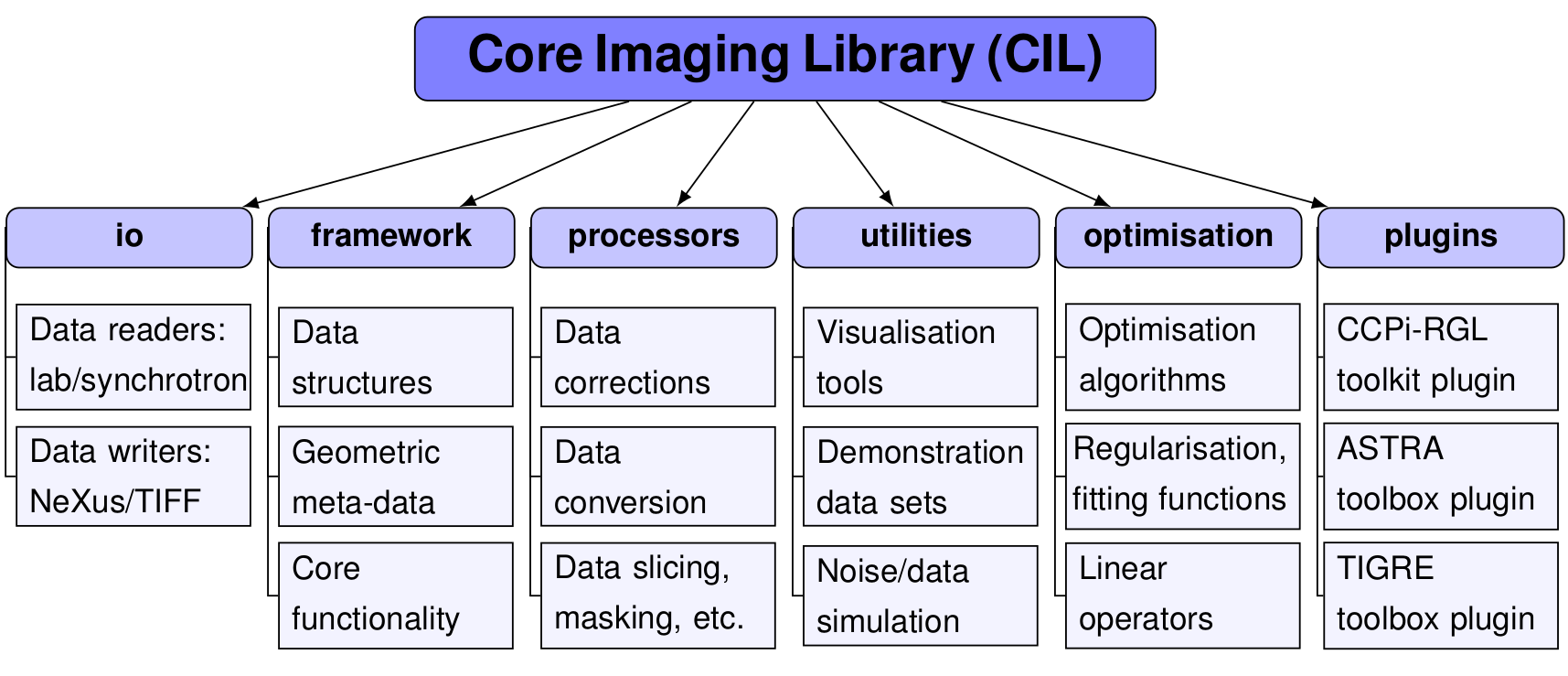}
\caption{Overview of CIL module structure and contents. The \textbf{cil.plugins} module contains wrapper code for other software and third-party libraries that need to be installed separately to be used by CIL. \label{fig:CILmodules}}
\end{figure}

 CIL is developed mainly in Python and binary distribution is currently via Anaconda. Instructions for installation and getting started are available at \url{https://www.ccpi.ac.uk/CIL} as well as at \cite{CILReleases}. The present version \cilversion{} consists of six modules, as shown in \cref{fig:CILmodules}.
 CIL is open-source software released under the Apache 2.0 license, while individual plugins may have a different license, e.g. \module{ccpi.plugins.astra} is GPLv3.
 In the following subsections the key functionality of each CIL module is explained and demonstrated, apart from \textbf{ccpi.optimisation} which is covered in \cref{sec:optimisation}.
 
As a running example (\cref{fig:DLSdata}) we employ a 3D parallel-beam X-ray CT data set from Beamline I13-2, Diamond Light Source, Harwell, UK. The sample consisted of a 0.5 mm aluminium cylinder with a piece of steel wire embedded in a small drilled hole. A droplet of salt water was placed on top, causing corrosion to form hydrogen bubbles. The data set, which was part of a fast time-lapse experiment, consists of 91 projections over 180$^\circ$, originally acquired as size 2560-by-2160 pixels, but provided in \cite{DLSdata} downsampled to 160-by-135 pixels.

\subsection{Data readers and writers} \label{sec:io}
Tomographic data comes in a variety of different formats depending on the instrument manufacturer or imaging facility. CIL currently supplies a native reader for Nikon's XTek data format, Zeiss' TXRM format, the NeXus format \cite{Konnecke2015} if exported by CIL, as well as TIFF stacks. Here ``native'' means that a CIL \code{AcquisitionData} object incl. \code{geometry} (as described in the following subsection) will be created by the CIL reader. Other data formats can be read using e.g. DXchange \cite{DeCarlo2014} and a CIL \code{AcquisitionData} object can be manually constructed.
CIL currently provides functionality to export/write data to disk in NeXus format or as a TIFF stack. 

The steel-wire dataset is included as an example in CIL. It is in NeXus format and can be loaded using \code{NEXUSDataReader}. For example data sets in CIL we provide a convenience method that saves the user from typing the path to the datafile:
\begin{center}
\begin{tcolorbox}[
    enhanced,
    attach boxed title to top center={yshift=-2mm},
    colback=darkspringgreen!20,
    colframe=darkspringgreen,
    colbacktitle=darkspringgreen,
    title=Load steel-wire example dataset,
    text width = 15cm,
    fonttitle=\bfseries\color{white},
    boxed title style={size=small,colframe=darkspringgreen,sharp corners},
    sharp corners,
]
\begin{minted}{python}
from cil.utilities.dataexample import SYNCHROTRON_PARALLEL_BEAM_DATA
data = SYNCHROTRON_PARALLEL_BEAM_DATA.get()
\end{minted}
\end{tcolorbox}
\end{center}

\subsection{Data structures, geometry and core functionality}

CIL provides two essential classes for data representation, namely  \code{AcquisitionData} for tomographic data and \code{ImageData} for reconstructed (or simulated) volume data. The steel-wire dataset was read in as an \code{AcquisitionData} that we can inspect with:
\begin{center}
\begin{tcolorbox}[
    enhanced,
    attach boxed title to top center={yshift=-2mm},
    colback=darkspringgreen!20,
    colframe=darkspringgreen,
    colbacktitle=darkspringgreen,
    title=Print AcquisitionData object to view essential properties,
    text width = 15cm,
    fonttitle=\bfseries\color{white},
    boxed title style={size=small,colframe=darkspringgreen,sharp corners},
    sharp corners,
]
\begin{minted}{python}
>>> print(data)
Number of dimensions: 3 
Shape: (91, 135, 160) 
Axis labels: ('angle', 'vertical', 'horizontal')
\end{minted}
\end{tcolorbox}
\end{center}

\begin{figure}[t]
\begin{minipage}{0.7\linewidth}
\centering
\includegraphics[height=4.35cm,clip,trim=0.55cm 0.1cm 2.1cm 1.0cm]{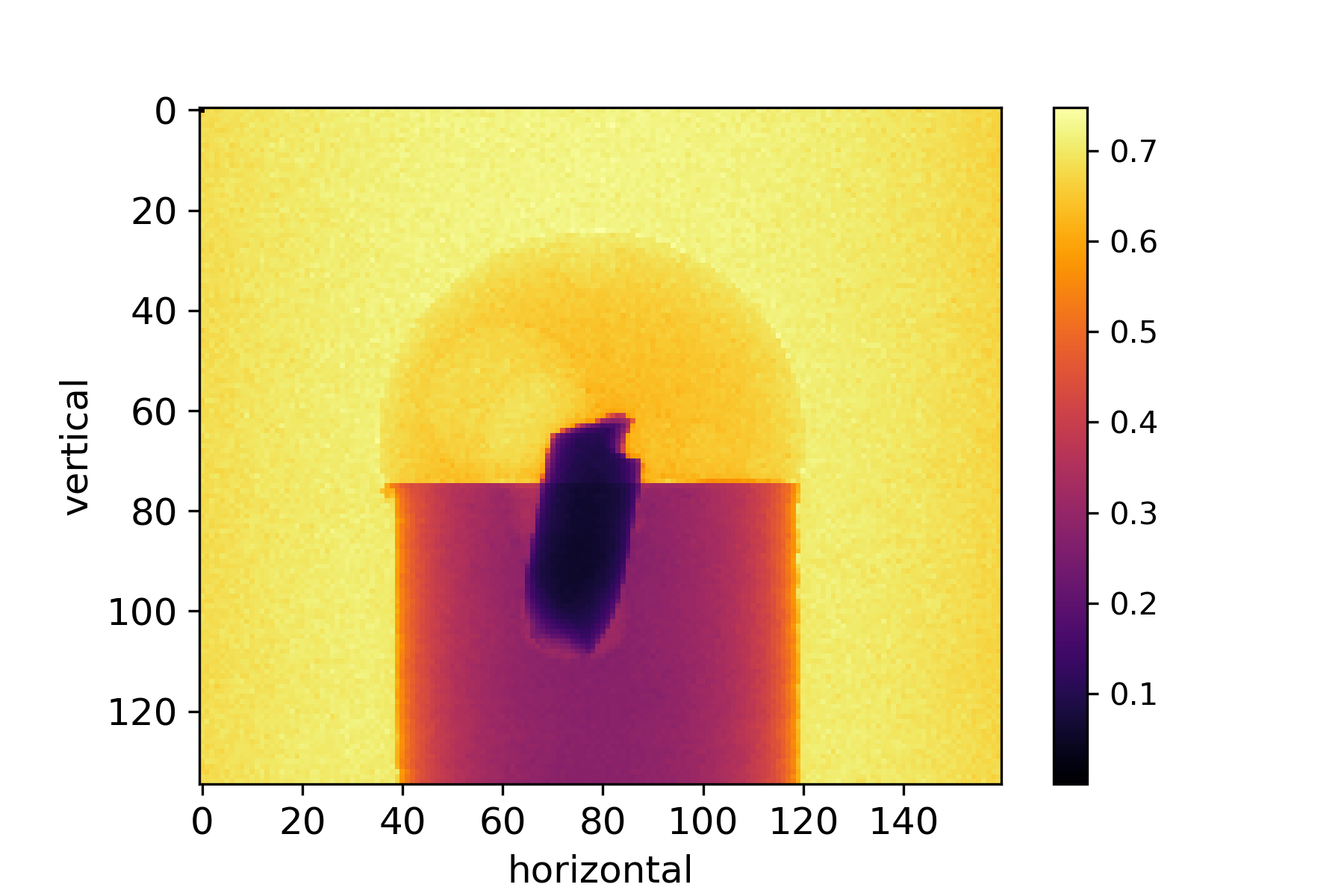}
\includegraphics[height=4.35cm,clip,trim=2.75cm 0.10cm 2.1cm 1.0cm]{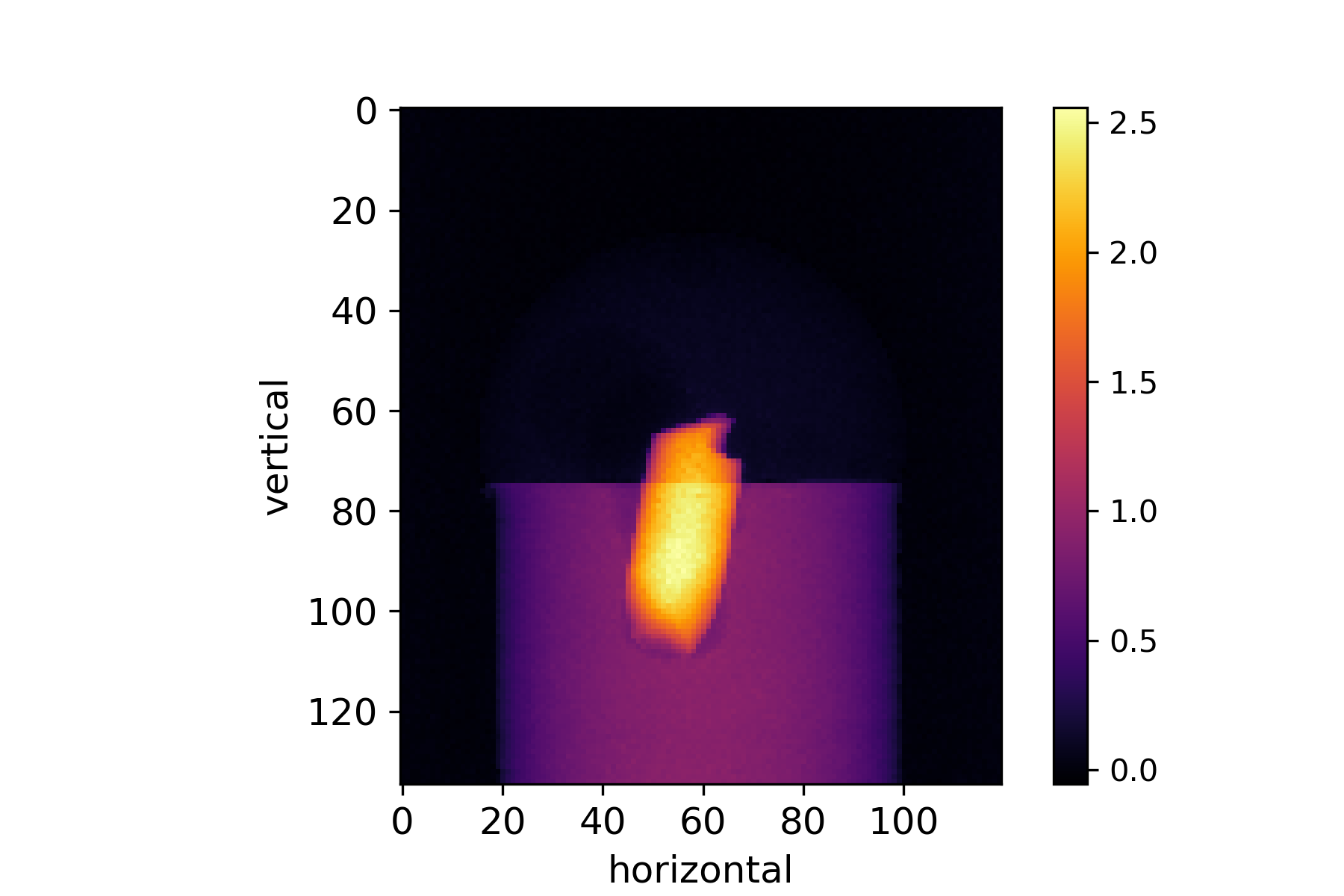}
\end{minipage}
\begin{minipage}{0.27\linewidth}
\centering
\includegraphics[width=\linewidth,clip,trim=1.2cm 0.4cm 2.25cm 0.75cm]{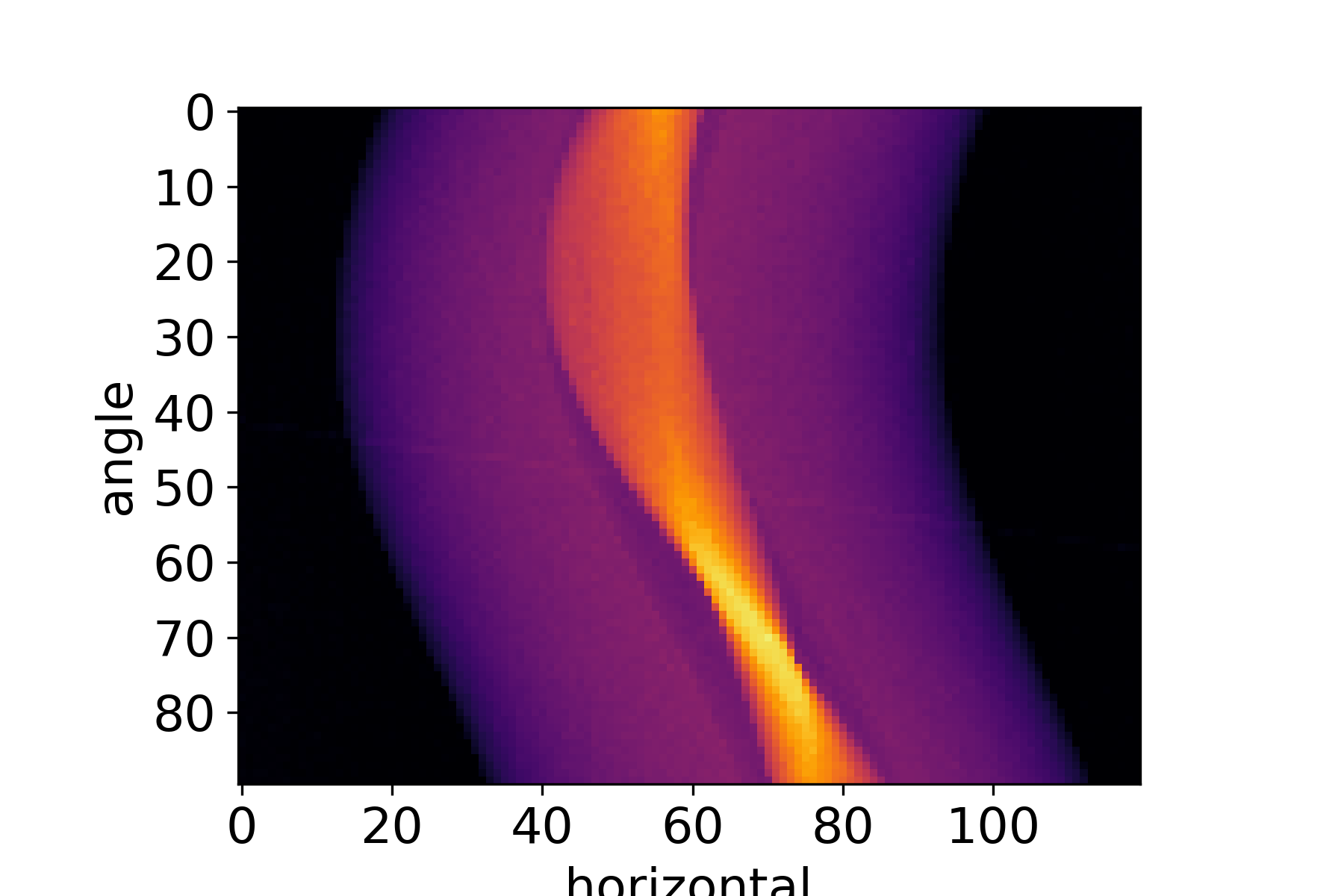}\\
\includegraphics[width=\linewidth,clip,trim=0.15cm 2.85cm 1.75cm 4.15cm]{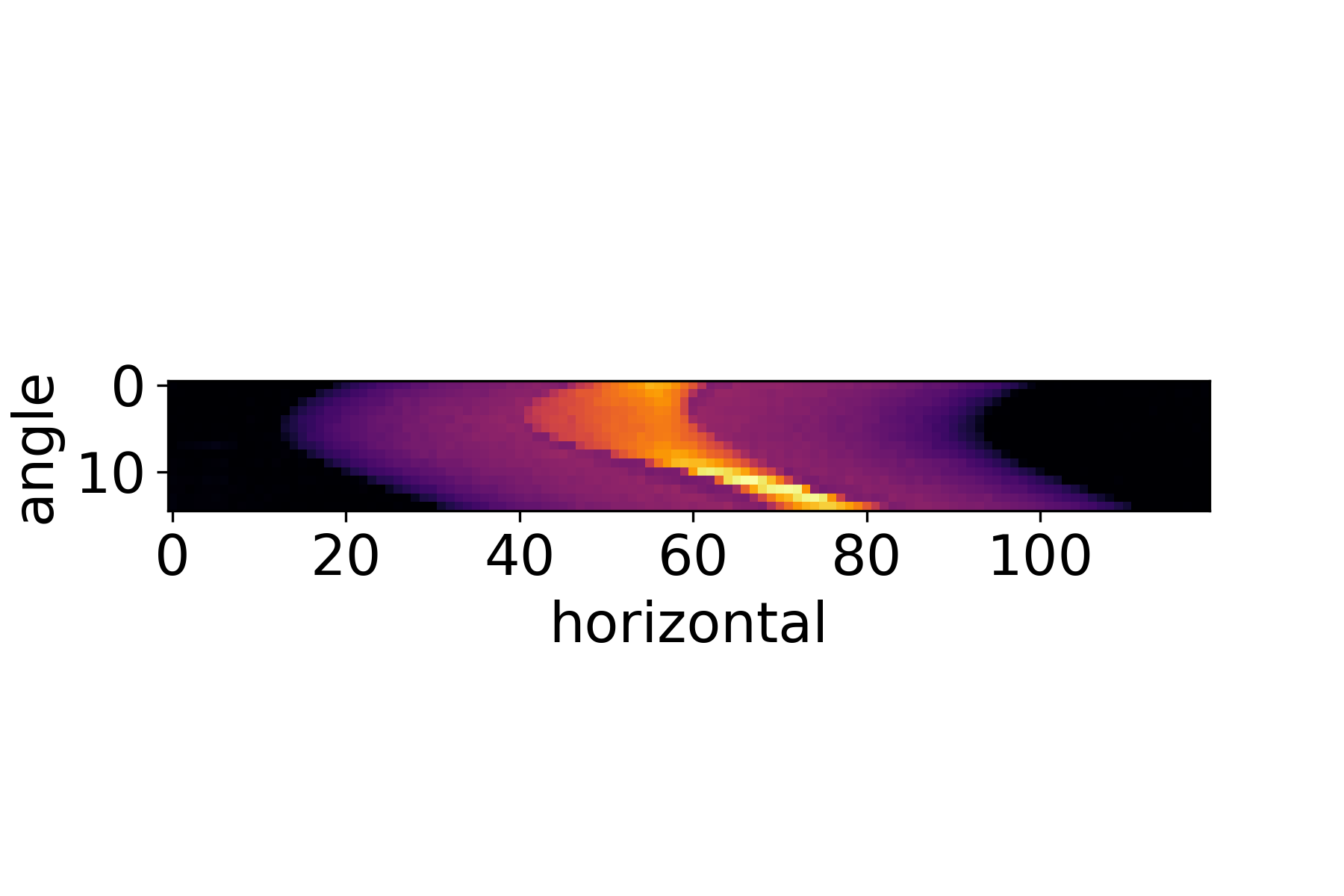}
\end{minipage}
\caption{Raw and preprocessed 3D parallel-beam X-ray CT steel-wire dataset. Left: Raw transmission projection. Centre: Scaled, cropped, centred and negative-log transformed projection. Right, top: Sinogram for slice \code{vertical=103}, all 90 angles. Right, bottom: Same, subsampled to 15 equi-spaced angles.}
\label{fig:DLSdata}
\end{figure}

At present, data is stored internally as a NumPy array and may be returned using the method \code{as_array()}. \code{AcquisitionData} and \code{ImageData} use string labels rather than a positional index to represent the dimensions. 
In the example data, \code{'angle'}, \code{'vertical'} and \code{'horizontal'} refer to 91 projections each with vertical size 135 and horizontal size 160. Labels enable the user to access subsets of data without knowing the details of how it is stored underneath. For example we can extract a single projection using the method \code{get_slice} with the label and display it (\cref{fig:DLSdata} left) as 
\begin{center}
\begin{tcolorbox}[
    enhanced,
    attach boxed title to top center={yshift=-2mm},
    colback=darkspringgreen!20,
    colframe=darkspringgreen,
    colbacktitle=darkspringgreen,
    title=Extract single projection and display as image,
    text width = 15cm,
    fonttitle=\bfseries\color{white},
    boxed title style={size=small,colframe=darkspringgreen,sharp corners},
    sharp corners,
]
\begin{minted}{python}
show2D(data.get_slice(angle=0), cmap='inferno', origin='upper-left')
\end{minted}
\end{tcolorbox}
\end{center}
where \code{show2D} is a display functiontter in \module{cil.utilities.display}. \code{show2D} displays dimension labels on plot axes as in \cref{fig:DLSdata}; subsequent plots omit these for space reasons.

Both \code{ImageData} and \code{AcquisitionData} behave much like a NumPy array with support for:
\begin{itemize}[noitemsep]
\item algebraic operators \code{+}, \code{-}, etc.,
\item relational operators \code{>}, \code{>=}, etc.,
\item common mathematical functions like \code{exp}, \code{log} and \code{abs}, \code{mean}, and
\item inner product \code{dot} and Euclidean norm \code{norm}.
\end{itemize}
This makes it easy to do a range of data processing tasks. For example in \cref{fig:DLSdata} (left) we note the projection (which is already flat-field normalised) has values around 0.7 in the background, and not 1.0 as in typical well-normalised data. This may lead to reconstruction artifacts. A quick-fix is to scale the image to have background value ca. 1.0. To do that we extract a row of the data toward the top, compute its mean and use it to normalise the data:
\begin{center}
\begin{tcolorbox}[
    enhanced,
    attach boxed title to top center={yshift=-2mm},
    colback=darkspringgreen!20,
    colframe=darkspringgreen,
    colbacktitle=darkspringgreen,
    title=Normalise data by mean over vertical slice of data,
    text width = 15cm,
    fonttitle=\bfseries\color{white},
    boxed title style={size=small,colframe=darkspringgreen,sharp corners},
    sharp corners,
]
\begin{minted}{python}
data = data / data.get_slice(vertical=20).mean()
\end{minted}
\end{tcolorbox}
\end{center}
Where possible in-place operations are supported to avoid unnecessary copying of data. For example the Lambert-Beer negative logarithm conversion can be done by: 
\begin{center}
\begin{tcolorbox}[
    enhanced,
    attach boxed title to top center={yshift=-2mm},
    colback=darkspringgreen!20,
    colframe=darkspringgreen,
    colbacktitle=darkspringgreen,
    title=In-place mathematical operations,
    text width = 15cm,
    fonttitle=\bfseries\color{white},
    boxed title style={size=small,colframe=darkspringgreen,sharp corners},
    sharp corners,
]
\begin{minted}{python}
data.log(out=data)
data *= -1
\end{minted}
\end{tcolorbox}
\end{center}

Geometric meta-data such as voxel dimensions and scan configuration is stored in \code{ImageGeometry} and \code{AcquisitionGeometry} objects available in the attribute \code{geometry} of \code{ImageData} and \code{AcquisitionData}. \code{AcquisitionGeometry} will normally be provided as part of an \code{AcquisitionData} produced by the CIL reader. It is also possible to manually create \code{AcquisitionGeometry} and \code{ImageGeometry} from a list of geometric parameters. Had the steel-wire dataset not had geometry information included, we could have set up its geometry with the following call:
\begin{center}
\begin{tcolorbox}[
    enhanced,
    attach boxed title to top center={yshift=-2mm},
    colback=darkspringgreen!20,
    colframe=darkspringgreen,
    colbacktitle=darkspringgreen,
    title=Manually define AcquisitionGeometry,
    text width = 15cm,
    fonttitle=\bfseries\color{white},
    boxed title style={size=small,colframe=darkspringgreen,sharp corners},
    sharp corners,
]
\begin{minted}{python}
ag = AcquisitionGeometry.create_Parallel3D()              \
        .set_panel(num_pixels=[160, 135])                 \
        .set_angles(angles=np.linspace(-88.2, 91.8, 91))
\end{minted}
\end{tcolorbox}
\end{center}
The first line creates a default 3D parallel-beam geometry with a rotation axis perpendicular to the beam propagation direction. The second and third lines specify the detector dimension and the angles at which projections are acquired. Numerous configuration options are available for bespoke geometries; this is illustrated in \cref{sec:lami}, see in particular \cref{fig:CILgeometrylami}, for an example of cone-beam laminography.
Similarly, \code{ImageGeometry} holds the geometric specification of a reconstructed volume, including numbers and sizes of voxels.

\subsection{Preprocessing data}

In CIL a \proc is a	 class that takes an \code{ImageData} or \code{AcquisitionData} as input, carries out some operations on it and returns an \code{ImageData} or \code{AcquisitionData}. Example uses include common preprocessing tasks such as resizing (e.g. cropping or binning/downsampling) data, flat-field normalization and correction for centre-of-rotation offset, see \cref{tab:processors} for an overview of \code{Processors} currently in CIL.

\begin{table}[t]
\centering
\caption{\proc{}s currently available in CIL. }
\label{tab:processors}
\begin{tabular}{ll}
\hline
\textbf{Name} & \textbf{Description}  \\
\hline
Binner & Downsample data in selected dimensions \\
CentreOfRotationCorrector & Find and correct for centre-of-rotation offset\\
Normaliser & Apply flat and dark field correction/normalisation \\
Padder & Pad/extend data in selected dimensions \\
Slicer & Extract data at specified indices\\
Masker & Apply binary mask to keep selected data only\\
MaskGenerator & Make binary mask to keep selected data only\\
RingRemover & Remove sinogram stripes to reduce ring artifacts \\
\hline
\end{tabular}
\vspace*{-4pt}
\end{table}

We will demonstrate centre-of-rotation correction and cropping using a \proc{}. Typically it is not possible to align the rotation axis perfectly with respect to the detector, and this leads to well-known centre-of-rotation reconstruction artifacts. CIL provides different techniques to estimate and compensate, the simplest being based on cross-correlation on the central slice. First the \proc instance must be created; this is an object instance which holds any parameters specified by the user; here which slice to operate on. Once created the \proc{} can carry out the processing task by calling it on the targeted data set. All this can be conveniently achieved in a single code line, as shown in the first line below.  

Afterwards, we use a \code{Slicer} to remove some of the empty parts of the projections by cropping 20 pixel columns on each side of all projections, while also discarding the final projection which is a mirror image of the first. This produces \code{data90}. We can further produce a subsampled data set \code{data15} by using another \code{Slicer}, keeping only every sixth projection.

\begin{center}
\begin{tcolorbox}[
    enhanced,
    attach boxed title to top center={yshift=-2mm},
    colback=darkspringgreen!20,
    colframe=darkspringgreen,
    colbacktitle=darkspringgreen,
    title=Use Processors to center and subsample data,
    text width = 15cm,
    fonttitle=\bfseries\color{white},
    boxed title style={size=small,colframe=darkspringgreen,sharp corners},
    sharp corners,
]
\begin{minted}{python}
data = CentreOfRotationCorrector.xcorr(slice_index='centre')(data)
data90 = Slicer(roi={'angle':(0,90), 'horizontal':(20,140)})(data)
data15 = Slicer(roi={'angle':(0,90,6)})(data90)
\end{minted}
\end{tcolorbox}
\end{center}

\Cref{fig:DLSdata} illustrates preprocessing and the final 90- and 15-projection sinograms; mainly the latter will be used in what follows to highlight differences between reconstruction methods.

\subsection{Auxiliary tools}
This module contains a number of useful tools:
\begin{itemize}
\item \module{dataexample}: Example data sets and test images such as the steel-wire dataset\footnote{Available from \url{https://github.com/TomographicImaging/CIL-Data}}.
\item \module{display}: Tools for displaying data as images, including the \code{show2D} used in the previous section and other interactive displaying tools for Jupyter notebooks.
\item \module{noise}: Tools to simulate different kinds of noise, including Gaussian and Poisson.
\item \module{quality\_measures}: Mathematical metrics Mean-Square-Error (MSE) and Peak-Signal-to-Noise-Ratio (PSNR) to quantify image quality against a ground-truth image.
\end{itemize}
Some of these tools are demonstrated in other sections of the present paper; for the rest we refer the reader to the CIL documentation.

\subsection{CIL Plugins and interoperability with SIRF}
CIL allows the use of third-party software through \emph{plugins} that wrap the desired functionality. At present the following three plugins are provided:
 \begin{itemize}
 \item \module{cil.plugins.ccpi\_regularisation} This plugin wraps a number of regularisation methods from the CCPi-RGL toolkit \cite{Kazantsev2019RGLTK} as CIL \fun{}s.
 \item \module{cil.plugins.astra}: This plugin provides access to CPU and GPU-accelerated forward and back projectors in ASTRA as well as the filtered back-projection (FBP) and Feldkamp-Davis-Kress (FDK) reconstruction methods for parallel and cone-beam geometries.
 \item \module{cil.plugins.tigre}: This plugin currently provides access to GPU-accelerated cone-beam forward and back projectors and the FDK reconstruction method of the TIGRE toolbox. 
 \end{itemize}
Furthermore, CIL is developed to be interoperable with the Synergistic Image Reconstruction Framework (SIRF) for PET and MR imaging \cite{SIRF2020}. This was achieved by synchronising naming conventions and basic class concepts:
\begin{itemize}
\item \module{sirf}: Data structures and acquisition models of SIRF can be used from CIL without a plugin, in particular with \module{cil.optimisation} one may specify and solve optimisation problems with SIRF data.  An example of this using PET data is given in \cref{sec:exsirf}.
\end{itemize}

\begin{figure}[t]
\includegraphics[height=4.0cm,clip,trim=2.35cm 0.0cm 3.6cm 1.0cm]{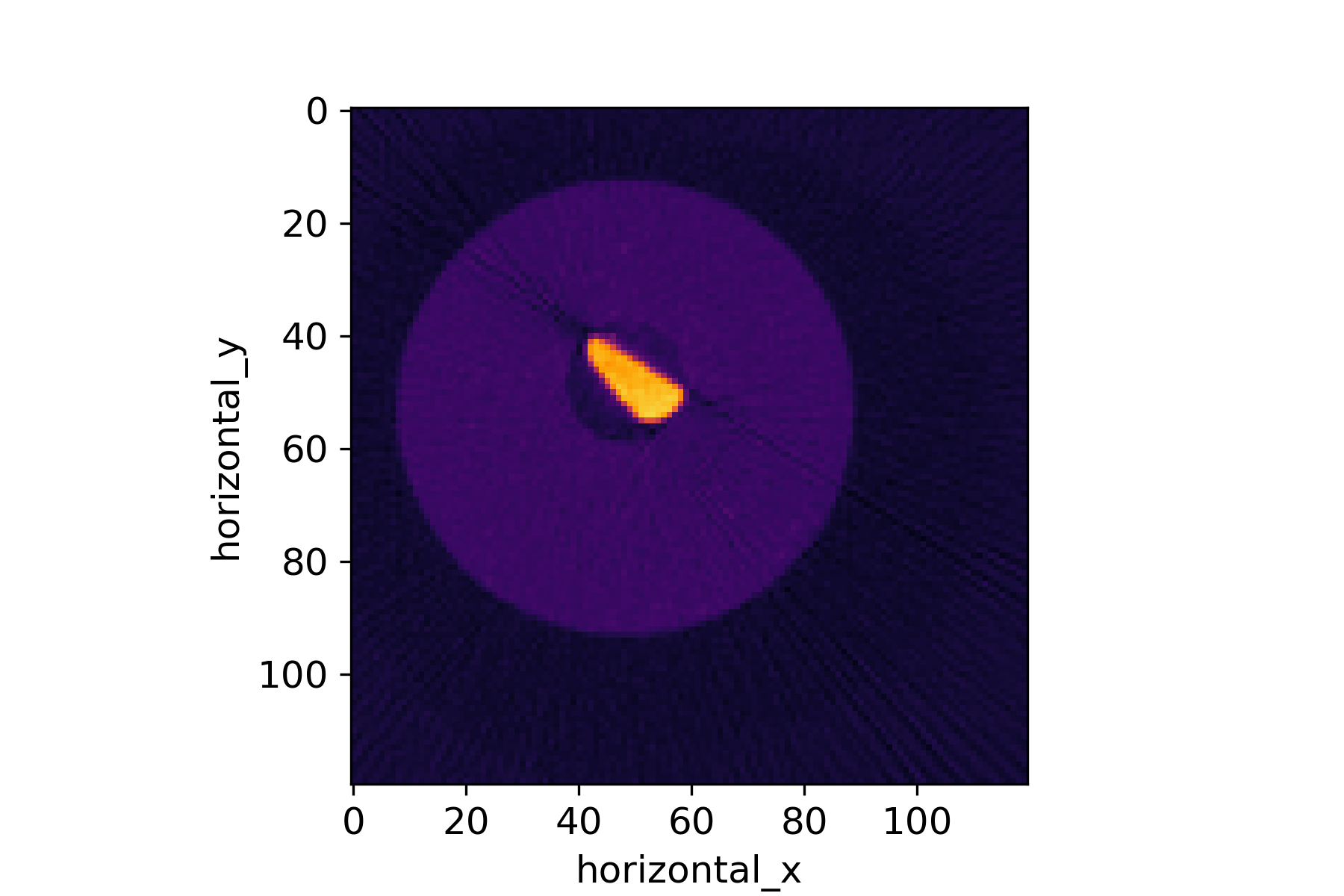}
\includegraphics[height=4.0cm,clip,trim=2.75cm 0.0cm 4cm 1.0cm]{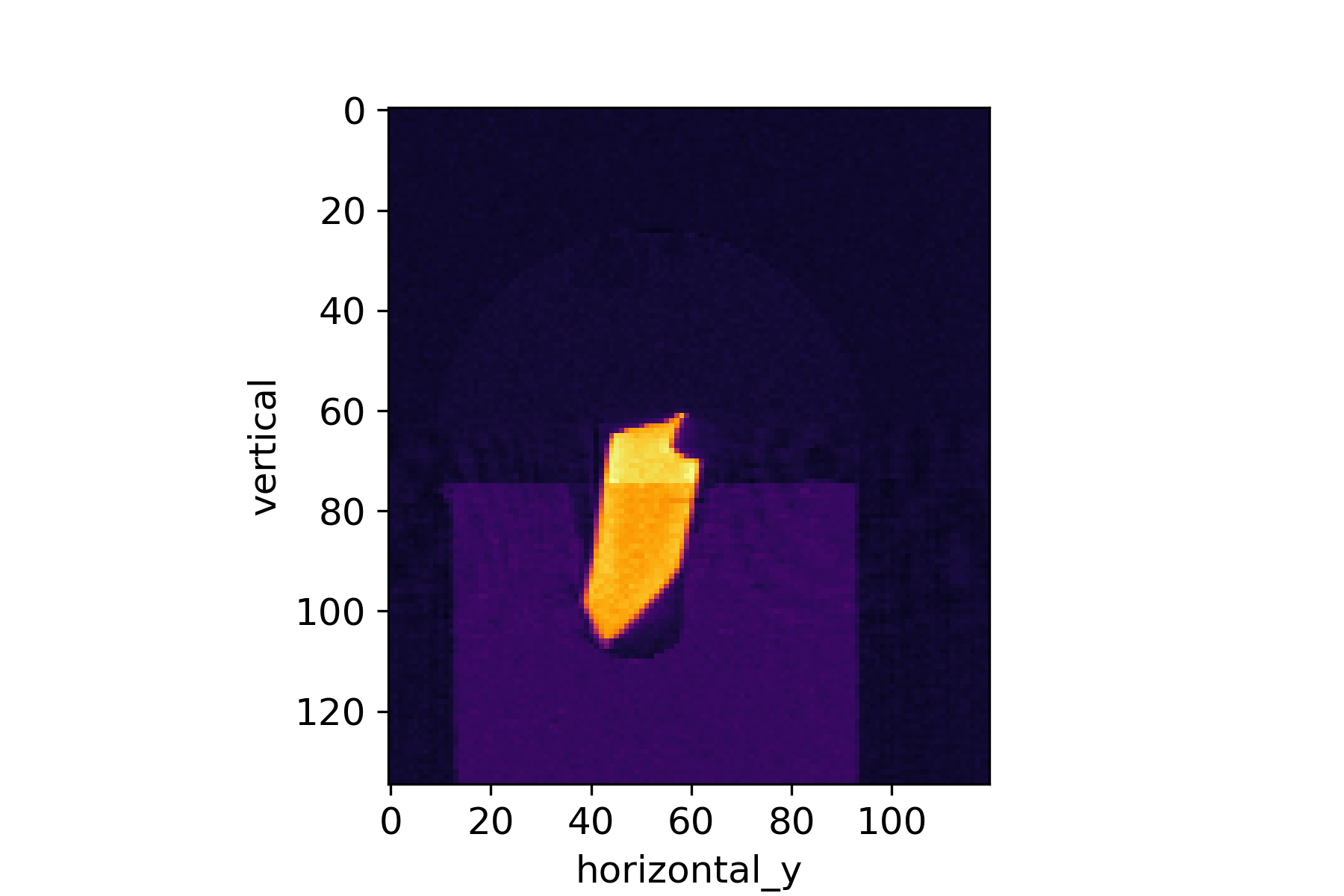}
\includegraphics[height=4.0cm,clip,trim=2.35cm 0.0cm 3.6cm 1.0cm]{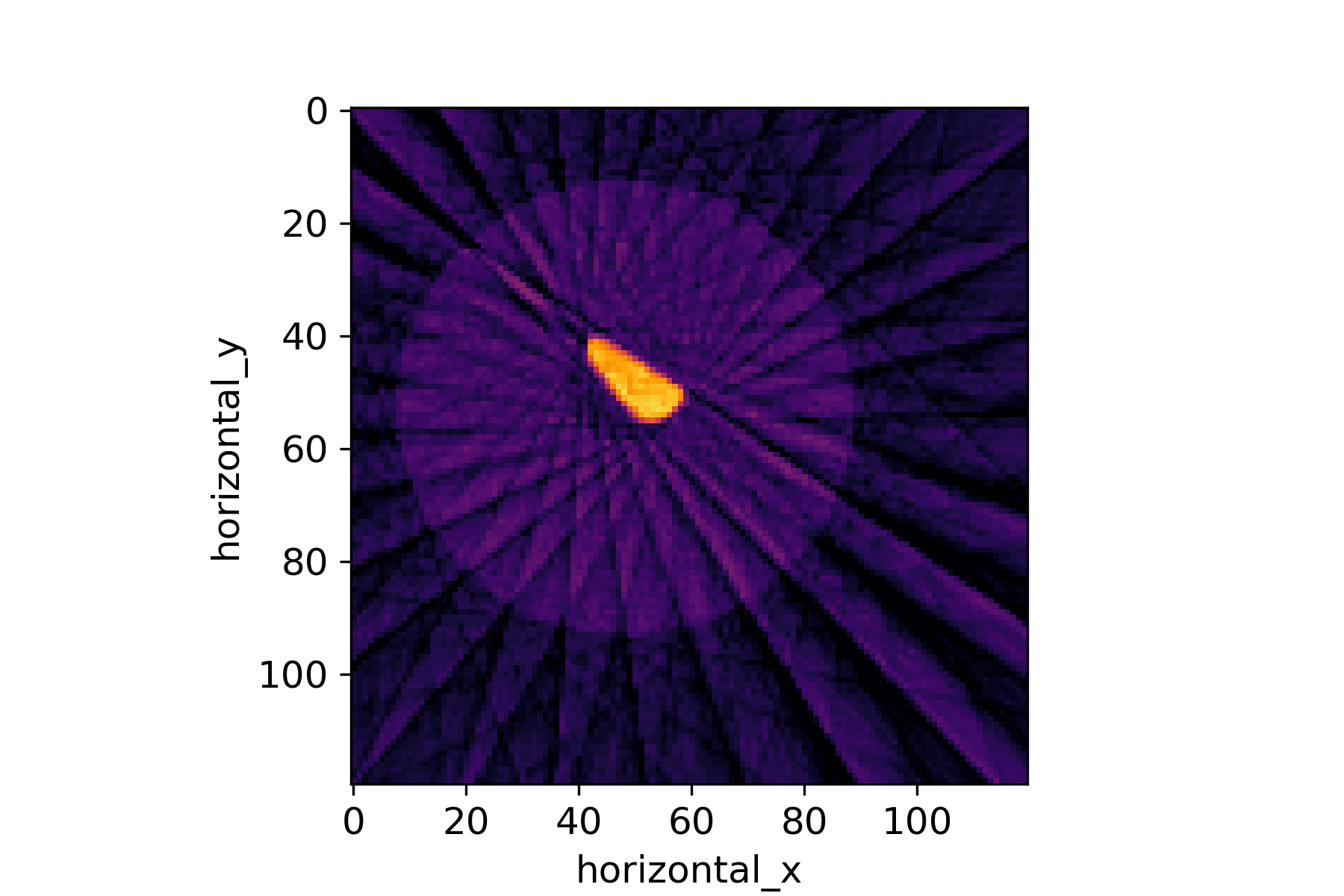}
\includegraphics[height=4.0cm,clip,trim=2.75cm 0.0cm 4cm 1.0cm]{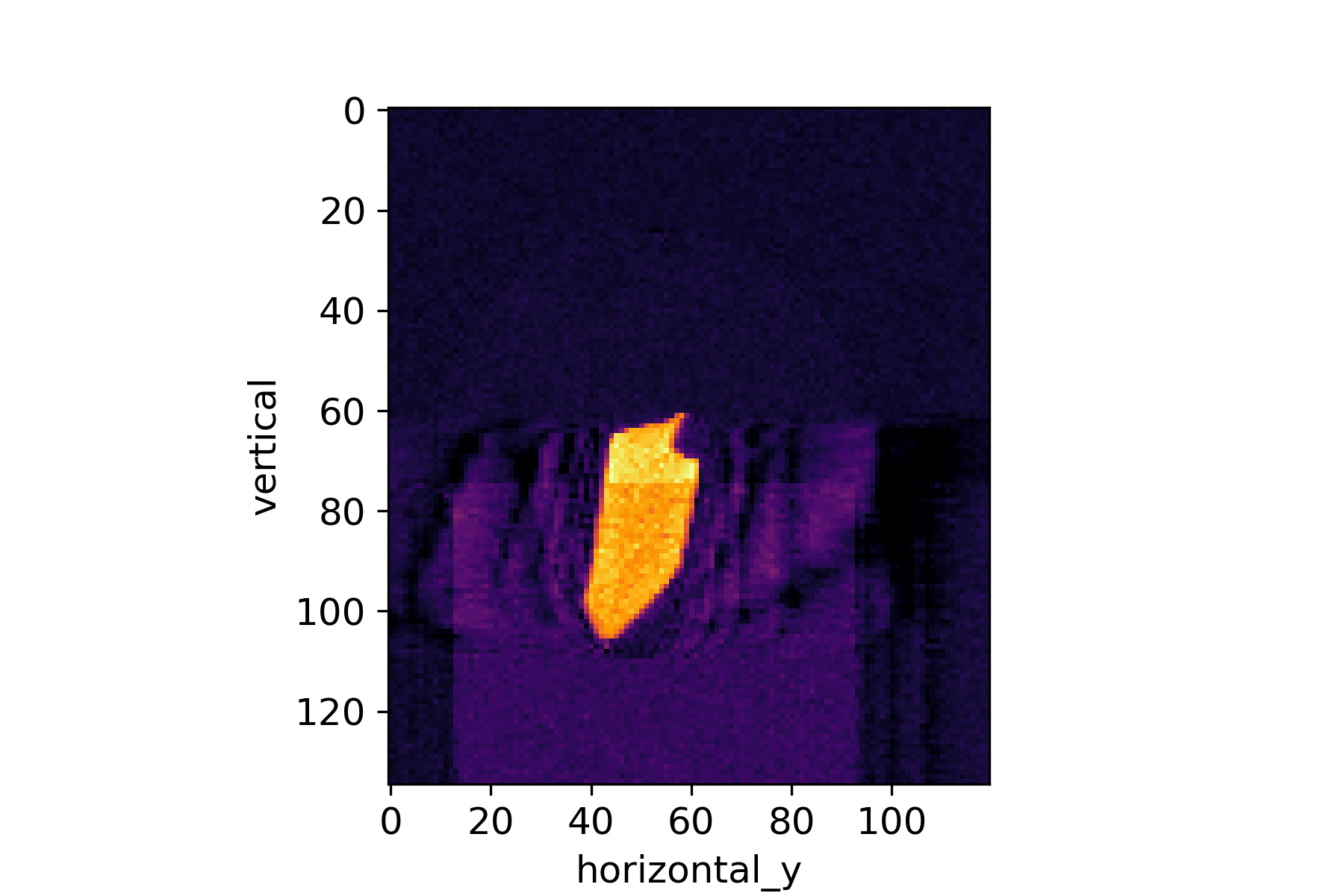}
\caption{Reconstructions of  steel-wire dataset by \code{FBP}. Left two: Horizontal and vertical slice using 90 projections. Right two: Same using 15 projections -- showing prominent streak artifacts. Colour range [-0.01, 0.11].}
\label{fig:DLSFBP}
\end{figure}

We demonstrate here how the \module{cil.plugins.astra} plugin, or \module{cil.plugins.tigre} plugin interchangeably, can be used to produce an FBP reconstruction of the steel-wire dataset using its \code{FBP} \proc{}.
To compute a reconstruction we must specify the geometry we want for the reconstruction volume; for convenience, a default \code{ImageGeometry} can be determined from a given \code{AcquisitionGeometry}. 
The \code{FBP} \proc{} can then be set up and in this instance we specify for it to use GPU-acceleration, and then call it on the data set to produce a reconstruction:
\begin{center}
\begin{tcolorbox}[
    enhanced,
    attach boxed title to top center={yshift=-2mm},
    colback=darkspringgreen!20,
    colframe=darkspringgreen,
    colbacktitle=darkspringgreen,
    title=Set up and run GPU-accelerated FBP algorithm from ASTRA plugin,
    text width = 15cm,
    fonttitle=\bfseries\color{white},
    boxed title style={size=small,colframe=darkspringgreen,sharp corners},
    sharp corners,
]
\begin{minted}{python}
data15.reorder(order='astra')
ag = data15.geometry
ig = ag.get_ImageGeometry()
recon = FBP(ig, ag, device='gpu')(data15)
\end{minted}
\end{tcolorbox}
\end{center}
The first line permutes the underlying data array to the specific dimension order required by \module{cil.plugins.astra}, which may differ from how data is read into CIL. 
Reconstructions for both the 90- and 15-projection steel-wire datasets are seen in \cref{fig:DLSFBP}, with notable streak artifacts in the subsampled case, as is typical with few projections. 

\section{Reconstruction by solving optimisation problems} \label{sec:optimisation}
FBP type reconstruction methods have very limited capability to model and address challenging data sets. For example the type and amount of noise cannot be modelled and prior knowledge such as non-negativity or smoothness cannot be incorporated. A much more flexible class of reconstruction methods arises from expressing the reconstructed image as the solution to an optimisation problem combining data and noise models and any prior knowledge.

The CIL optimisation module makes it simple to specify a variety of optimisation problems for reconstruction and provides a range of optimisation algorithms for their solution. 

\subsection{Operators}

The \textbf{ccpi.optimisation} module is built around the generic linear inverse problem
\begin{align} \label{eq:linearinverseproblem}
\sysmat \im = \sino,
\end{align}
where $\sysmat$ is a linear operator, $\im$ is the image to be determined, and $\sino$ is the measured data. In CIL $\im$ and $\sino$ are normally represented by \code{ImageData} and \code{AcquisitionData} respectively, and $\sysmat$ by a \code{LinearOperator}. The spaces that a \code{LinearOperator} maps from and to are represented in attributes \code{domain} and \texttt{range}; these should each hold an \code{ImageGeometry} or \code{AcquisitionGeometry} that match with that of $u$ and $b$, respectively.

Reconstruction methods rely on two essential methods of a \code{LinearOperator}, namely \code{direct}, which evaluates $Av$ for a given $v$, and \code{adjoint}, which evaluates $A^*z$ for a given $z$, where $A^*$ is the adjoint operator of $A$. For example, in a  \code{LinearOperator} representing the discretised Radon transform for tomographic imaging, \code{direct} is \emph{forward projection}, i.e., computing the sinogram corresponding to a given image, while \code{adjoint} corresponds to  \emph{back-projection}.

\Cref{tab:operators} provides an overview of the \op{}s available in the current version of CIL. It includes imaging models such as \code{BlurringOperator} for image deblurring problems and mathematical operators such as \code{IdentityOperator} and \code{GradientOperator} to act as building blocks for specifying optimisation problems. 
\op{}s can be combined to create new \op{}s through addition, scalar multiplication and composition.

The bottom two row contains \code{ProjectionOperator} from both \module{cil.plugins.astra} and \module{cil.plugins.tigre}, which wraps forward and back-projectors from the ASTRA and TIGRE toolboxes respectively, and can be used interchangeably. A \code{ProjectionOperator} can be set up simply by
\begin{center}
\begin{tcolorbox}[
    enhanced,
    attach boxed title to top center={yshift=-2mm},
    colback=darkspringgreen!20,
    colframe=darkspringgreen,
    colbacktitle=darkspringgreen,
    title=Create ProjectionOperator from image and acquisition geometries,
    text width = 15cm,
    fonttitle=\bfseries\color{white},
    boxed title style={size=small,colframe=darkspringgreen,sharp corners},
    sharp corners,
]
\begin{minted}{python}
A = ProjectionOperator(ig, ag)
\end{minted}
\end{tcolorbox}
\end{center}
and from the \code{AcquisitionGeometry} provided the relevant 2D or 3D, parallel-beam or cone-beam geometry employed; in case of the steel-wire dataset, a 3D parallel-beam geometry.

\begin{table}[tb]
\centering
\caption{\op{}s in CIL; and \op{}s from \module{cil.plugins.astra} and \module{cil.plugins.tigre} in bottom two rows.}
\label{tab:operators}
\begin{tabular}{lllll}
\hline
\textbf{Name}  & \textbf{Description} \\
\hline
BlockOperator & Form block (array) operator from multiple operators\\
BlurringOperator & Apply point spread function to blur an image\\
ChannelwiseOperator & Apply the same Operator to all channels \\
DiagonalOperator & Form a diagonal operator from image/acquisition data \\
FiniteDifferenceOperator & Apply finite differences in selected dimension \\
GradientOperator & Apply finite difference to multiple/all dimensions \\
IdentityOperator & Apply identity operator, i.e., return input  \\
MaskOperator & From binary input, keep selected entries, mask out rest\\
SymmetrisedGradientOperator & Apply symmetrised gradient, used in TGV \\
ZeroOperator & Operator of all zeroes \\
\hline
ProjectionOperator & Tomography forward/back-projection from ASTRA\\ 
ProjectionOperator & Tomography forward/back-projection from TIGRE\\ \hline
\end{tabular}
\vspace*{-4pt}
\end{table}



\subsection{Algebraic iterative reconstruction methods}

One of the most basic optimisation problems for reconstruction is least-squares minimisation,
\begin{equation}
\im^\star = \argmin_\im \|\sysmat \im - \sino\|_2^2,
\end{equation}
where we seek to find the image $\im$ that fits the data the best, i.e., in which the norm of the residual $\sysmat \im - \sino$ takes on the smallest possible value; this $\im$ we denote $\im^\star$ and take as our reconstruction.

The Conjugate Gradient Least Squares (CGLS) algorithm \cite{CGLS} is an algebraic iterative method that solves exactly this problem. In CIL it is available as \code{CGLS}, which is an example of an \code{Algorithm} object. The following code sets up a \code{CGLS} algorithm -- inputs required are an initial image, the operator (here \code{ProjectionOperator} from \module{cil.plugins.astra}), the data and an upper limit on the number iterations to run -- and runs a specified number of iterations with verbose printing:
\begin{center}
\begin{tcolorbox}[
    enhanced,
    attach boxed title to top center={yshift=-2mm},
    colback=darkspringgreen!20,
    colframe=darkspringgreen,
    colbacktitle=darkspringgreen,
    title=Set up and run CGLS algorithm,
    text width = 15cm,
    fonttitle=\bfseries\color{white},
    boxed title style={size=small,colframe=darkspringgreen,sharp corners},
    sharp corners,
]
\begin{minted}{python}
x0 = ig.allocate(0.0)
b = data15
myCGLS = CGLS(initial=x0, operator=A, data=b, max_iteration=1000)
myCGLS.run(20, verbose=1)
\end{minted}
\end{tcolorbox}
\end{center}
At this point the reconstruction is available as \code{myCGLS.solution} and can be displayed or otherwise analysed. The object-oriented design of \code{Algorithm} means that iterating can be resumed from the current state, simply by another \code{myCGLS.run} call.

As imaging operators are often ill-conditioned with respect to inversion, small errors and inconsistencies tend to magnify during  the solution process, typically rendering the final least squares $\im^\star$ useless. CGLS exhibits semi-convergence \cite{PCHDIP} meaning that in the initial iterations the solution will approach the true underlying solution, but from a certain point the noise will increasingly contaminate the solution. The number of iterations therefore has an important regularising effect and must be chosen with care. 

\begin{figure}[t]
\includegraphics[height=4.1cm,clip,trim=3.7cm 1.3cm 3.9cm 1.25cm]{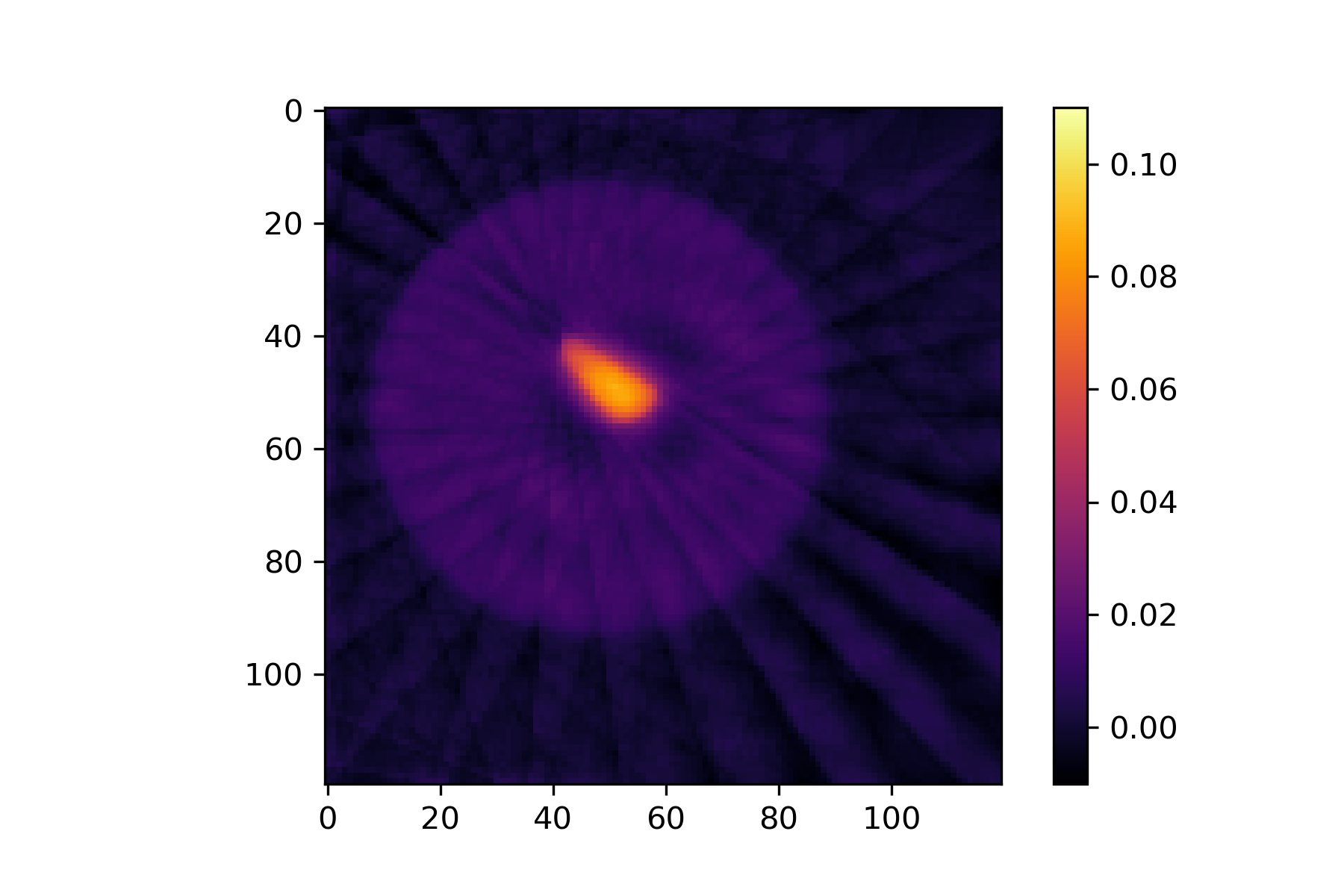}
\includegraphics[height=4.1cm,clip,trim=4.55cm 1.3cm 3.8cm 1.25cm]{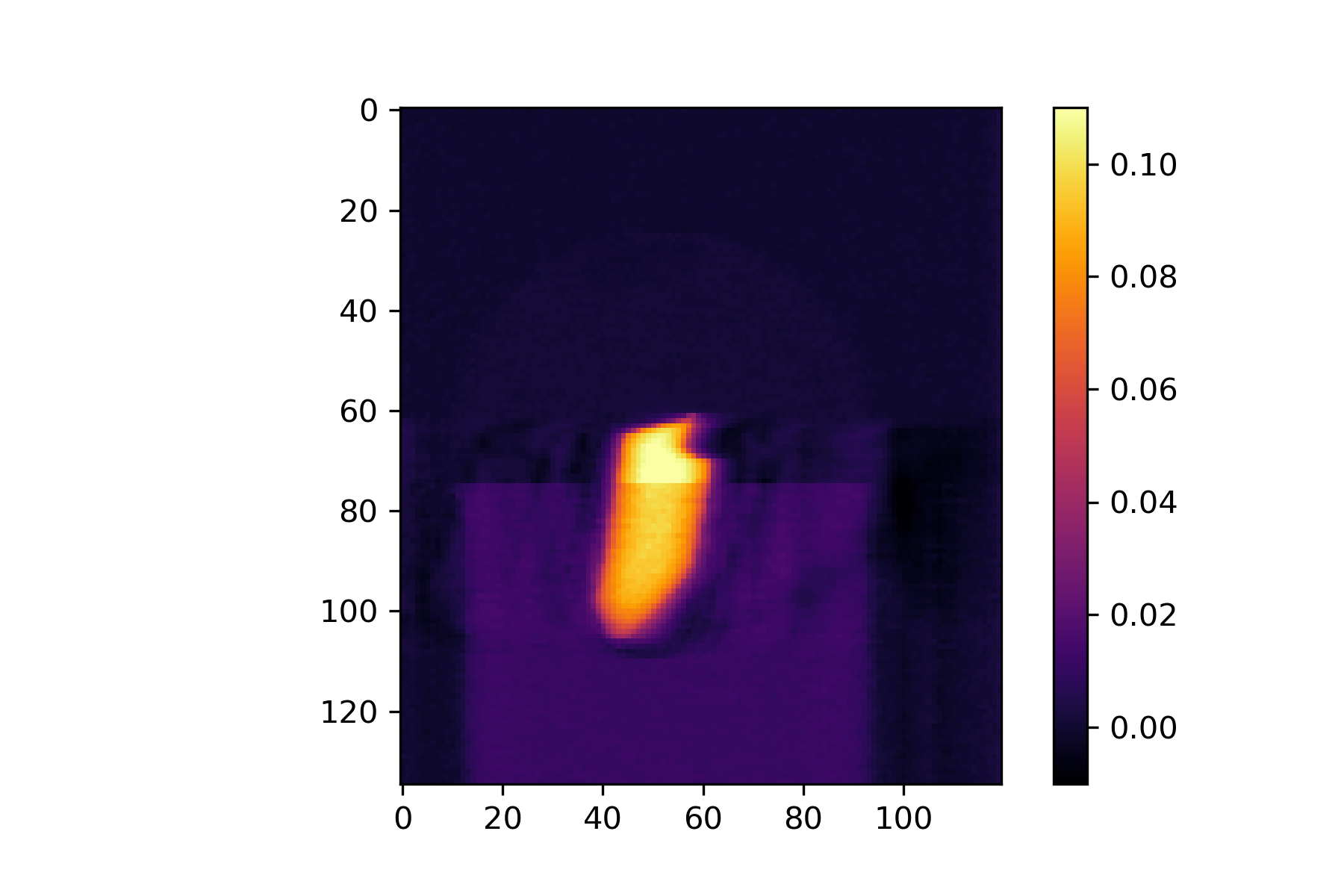}
\includegraphics[height=4.1cm,clip,trim=3.7cm 1.3cm 3.9cm 1.25cm]{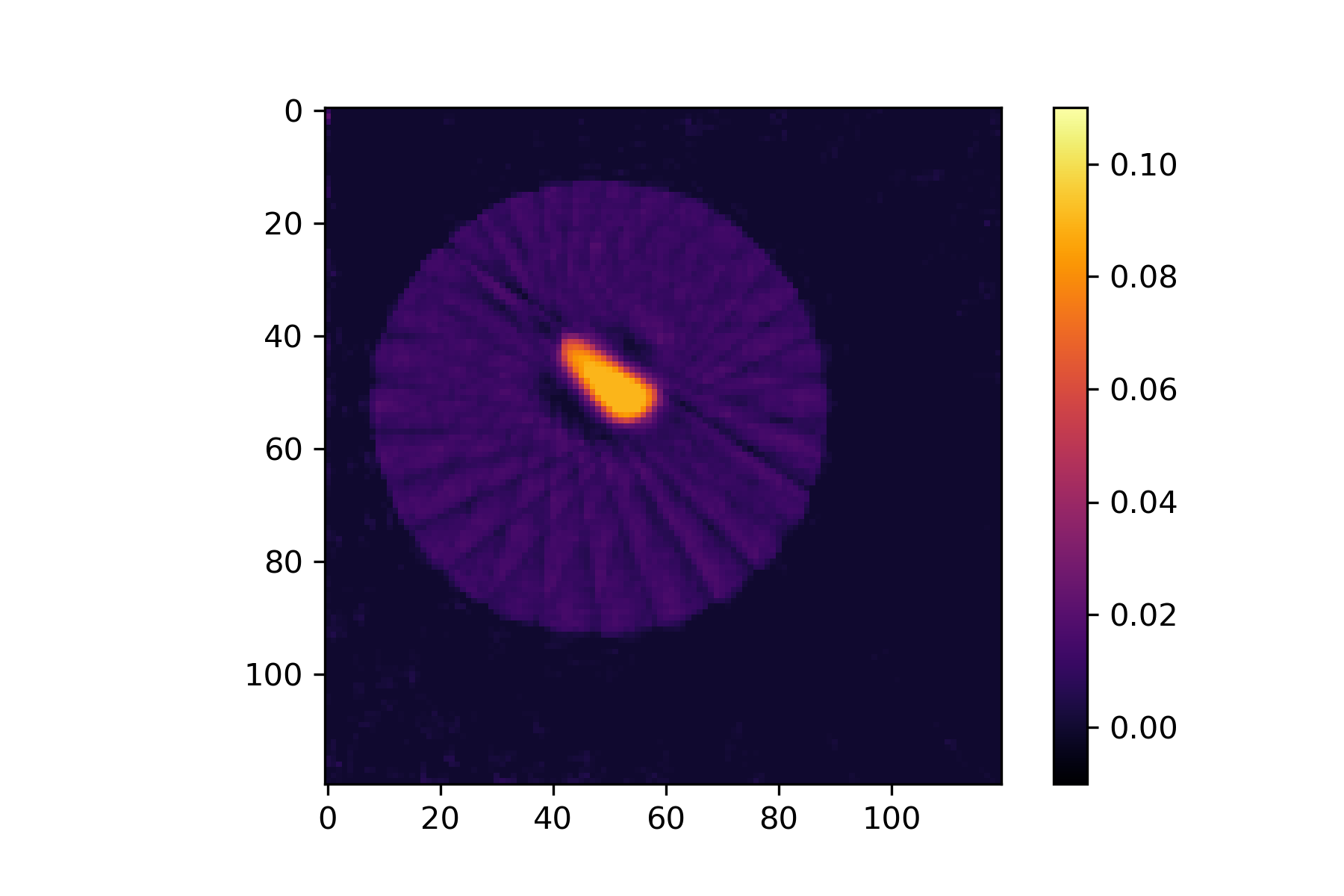}
\includegraphics[height=4.1cm,clip,trim=4.55cm 1.3cm 3.8cm 1.25cm]{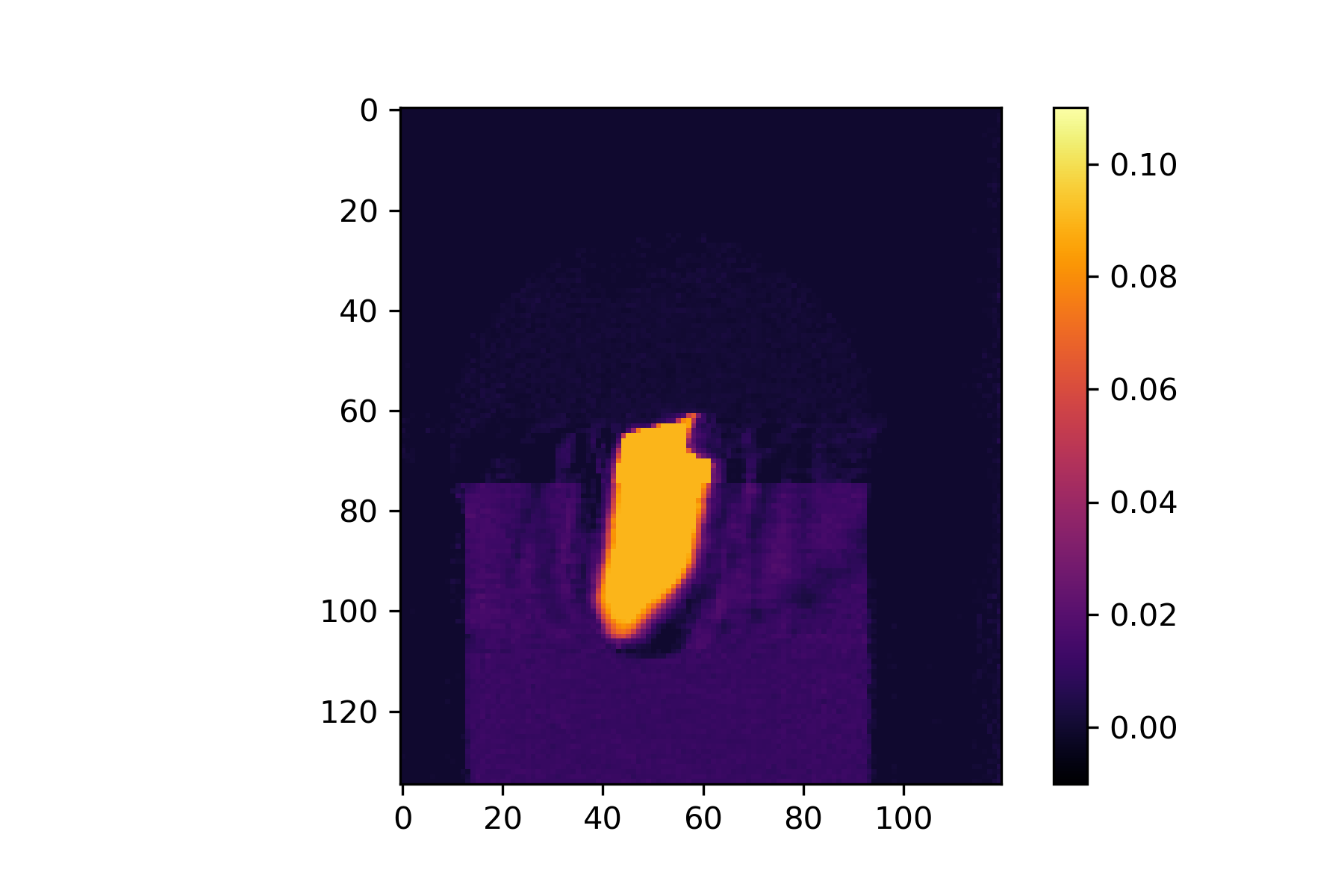}
\caption{Algebraic iterative reconstruction of 15-projection 3D steel-wire dataset. Left two: Horizontal and vertical slices, 20-iteration \code{CGLS} reconstruction. Right two: Same using \code{SIRT}, lower/upper bounds 0.0/0.09. Colour range [-0.01,0.11].}
\label{fig:DLScglssirt}
\end{figure}

CIL also provides the Simultaneous Iterative Reconstruction Technique (SIRT) as \code{SIRT}, which solves a particular weighted least-squares problem \cite{GILBERT1972,airtoolsii2018}. As with CGLS, it exhibits semi-convergence, however tends to require more iterations. An advantage of SIRT is that it admits the specification of convex constraints, such as a box constraints (upper and lower bounds) on $u$; this is done using optional input arguments \code{lower} and \code{upper}:
\begin{center}
\begin{tcolorbox}[
    enhanced,
    attach boxed title to top center={yshift=-2mm},
    colback=darkspringgreen!20,
    colframe=darkspringgreen,
    colbacktitle=darkspringgreen,
    title=Set up and run SIRT algorithm with bounds on pixel values,
    text width = 15cm,
    fonttitle=\bfseries\color{white},
    boxed title style={size=small,colframe=darkspringgreen,sharp corners},
    sharp corners,
]
\begin{minted}{python}
mySIRT = SIRT(initial=x0, operator=A, data=b, max_iteration=1000, \
              lower=0.0, upper=0.09)
mySIRT.run(200, verbose=1)
\end{minted}
\end{tcolorbox}
\end{center}
In \cref{fig:DLScglssirt} we see that \code{CGLS} reduces streaks but blurs edges. \code{SIRT} further reduces streaks and sharpens edges to the background; this is an effect of the nonnegativity constraint. In the steel wire example data the upper bound of 0.09 is attained causing a more uniform appearance with sharper edges.

\subsection{Tikhonov regularisation with BlockOperator and BlockDataContainer} \label{sec:tikhonovblock}
Algebraic iterative methods like CGLS and SIRT enforce regularisation of the solution implicitly by terminating iterations early. A more explicit form of regularisation is to include it directly in an optimisation formulation. The archetypal such method is Tikhonov regularisation which takes the form
\begin{equation} \label{eq:tikhonov}
\im^\star = \argmin_\im \left\{ \|\sysmat \im - \sino \|_2^2 + \alpha^2 \|D\im\|_2^2 \right\},
\end{equation}
where $D$ is some operator, the properties of which govern the appearance of the solution. In the simplest form $D$ can be taken as the identity operator. Another common choice is a discrete gradient implemented as a finite-difference operator.
The \emph{regularisation parameter} $\alpha$ governs the balance between the data fidelity term and the regularisation term. Conveniently, Tikhonov regularisation can be analytically rewritten as an equivalent least-squares problem, namely
\begin{equation} \label{eq:tikhonovls}
\im^\star = \argmin_\im  \left\|\tilde{\sysmat} \im - \tilde{\sino} \right\|_2^2, \qquad \text{where} \qquad \tilde{\sysmat} = \binom{\sysmat}{\alpha D} \quad \text{and} \quad \tilde{\sino} =  \binom{\sino}{0},
\end{equation}
where the $0$ corresponds to the range of $D$. We can use the CGLS algorithm to solve \cref{eq:tikhonovls} but we need a way to express the block structure of $\tilde{\sysmat}$ and $\tilde{\sino}$. This is achieved by the \code{BlockOperator} and \code{BlockDataContainer} of CIL:
\begin{center}
\begin{tcolorbox}[
    enhanced,
    attach boxed title to top center={yshift=-2mm},
    colback=darkspringgreen!20,
    colframe=darkspringgreen,
    colbacktitle=darkspringgreen,
    title=Set up Tikhonov regularisation for CGLS using BlockOperator and BlockDataContainer,
    text width = 15cm,
    fonttitle=\bfseries\color{white},
    boxed title style={size=small,colframe=darkspringgreen,sharp corners},
    sharp corners,
]
\begin{minted}{python}
alpha = 1.0
D = IdentityOperator(ig)
Atilde = BlockOperator(A, alpha*D)
z = D.range.allocate(0.0)
btilde = BlockDataContainer(b, z)
\end{minted}
\end{tcolorbox}
\end{center}
 If instead we want the discrete gradient as $D$ we simply replace the second line by:
\begin{center}
\begin{tcolorbox}[
    enhanced,
    attach boxed title to top center={yshift=-2mm},
    colback=darkspringgreen!20,
    colframe=darkspringgreen,
    colbacktitle=darkspringgreen,
    title=Set up GradientOperator for use in regularisation,
    text width = 15cm,
    fonttitle=\bfseries\color{white},
    boxed title style={size=small,colframe=darkspringgreen,sharp corners},
    sharp corners,
]
\begin{minted}{python}
D = GradientOperator(ig)
\end{minted}
\end{tcolorbox}
\end{center}
\code{GradientOperator} automatically works out from the \code{ImageGeometry} \code{ig} which dimensions are available and sets up finite differencing in all dimensions. If two or more dimensions are present, \code{D} will in fact be a \code{BlockOperator} with a finite-differencing block for each dimension. CIL supports nesting of a \code{BlockOperator} inside another, so that Tikhonov regularisation with a \code{Gradient} operator can be conveniently expressed.

\begin{figure}[t]
\includegraphics[height=4.1cm,clip,trim=3.7cm 1.3cm 3.9cm 1.25cm]{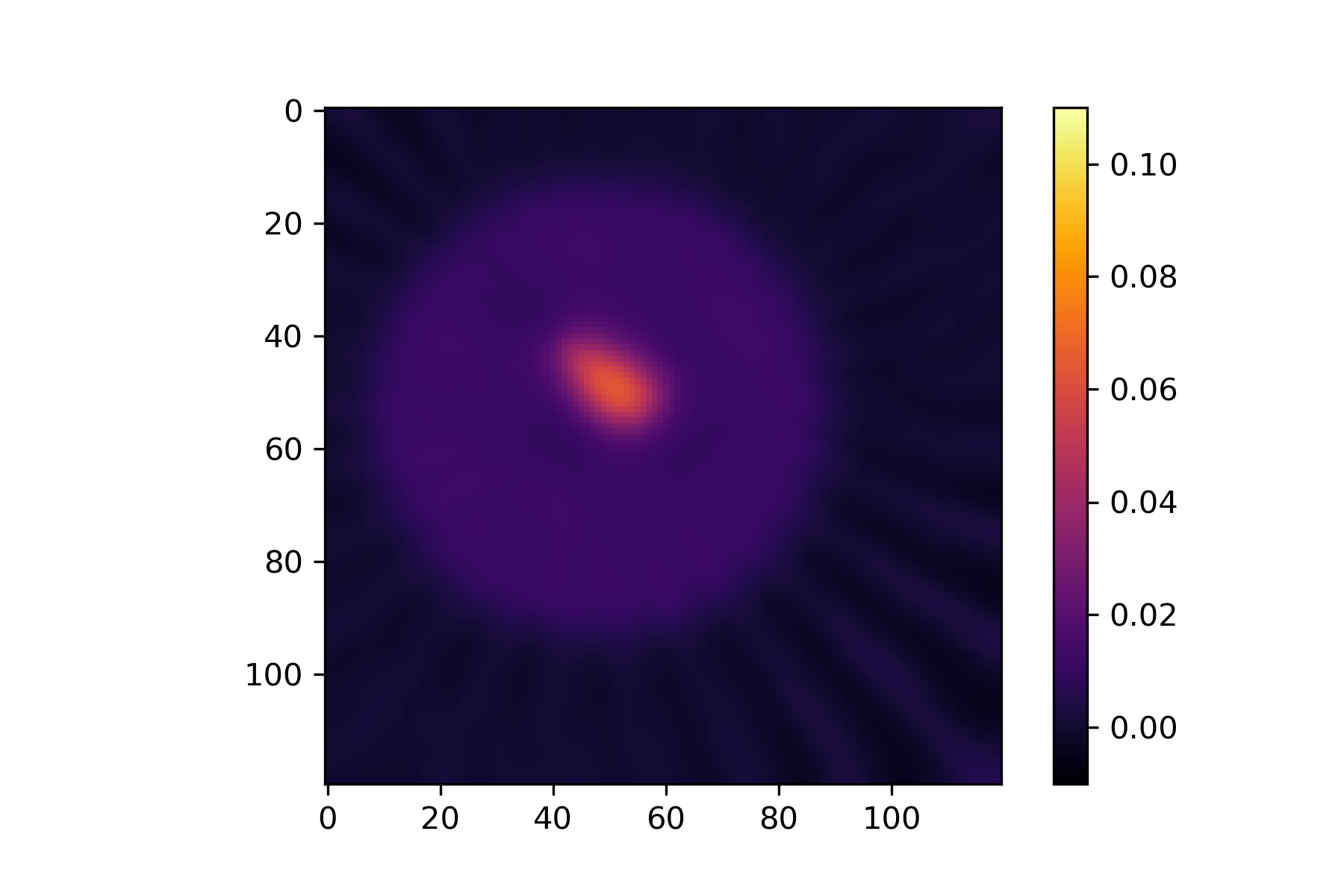}
\includegraphics[height=4.1cm,clip,trim=4.55cm 1.3cm 3.8cm 1.25cm]{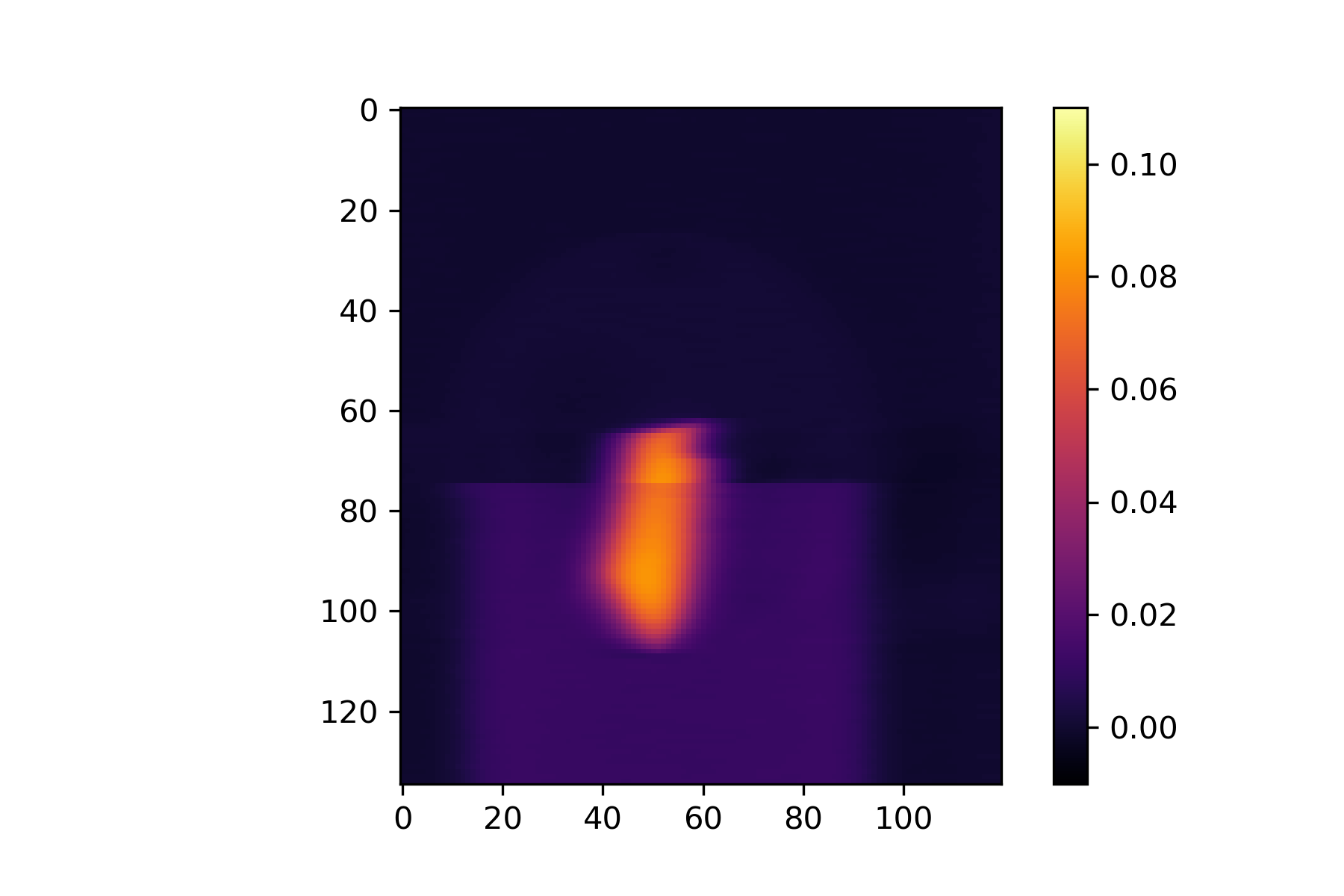}
\includegraphics[height=4.1cm,clip,trim=3.7cm 1.3cm 3.9cm 1.25cm]{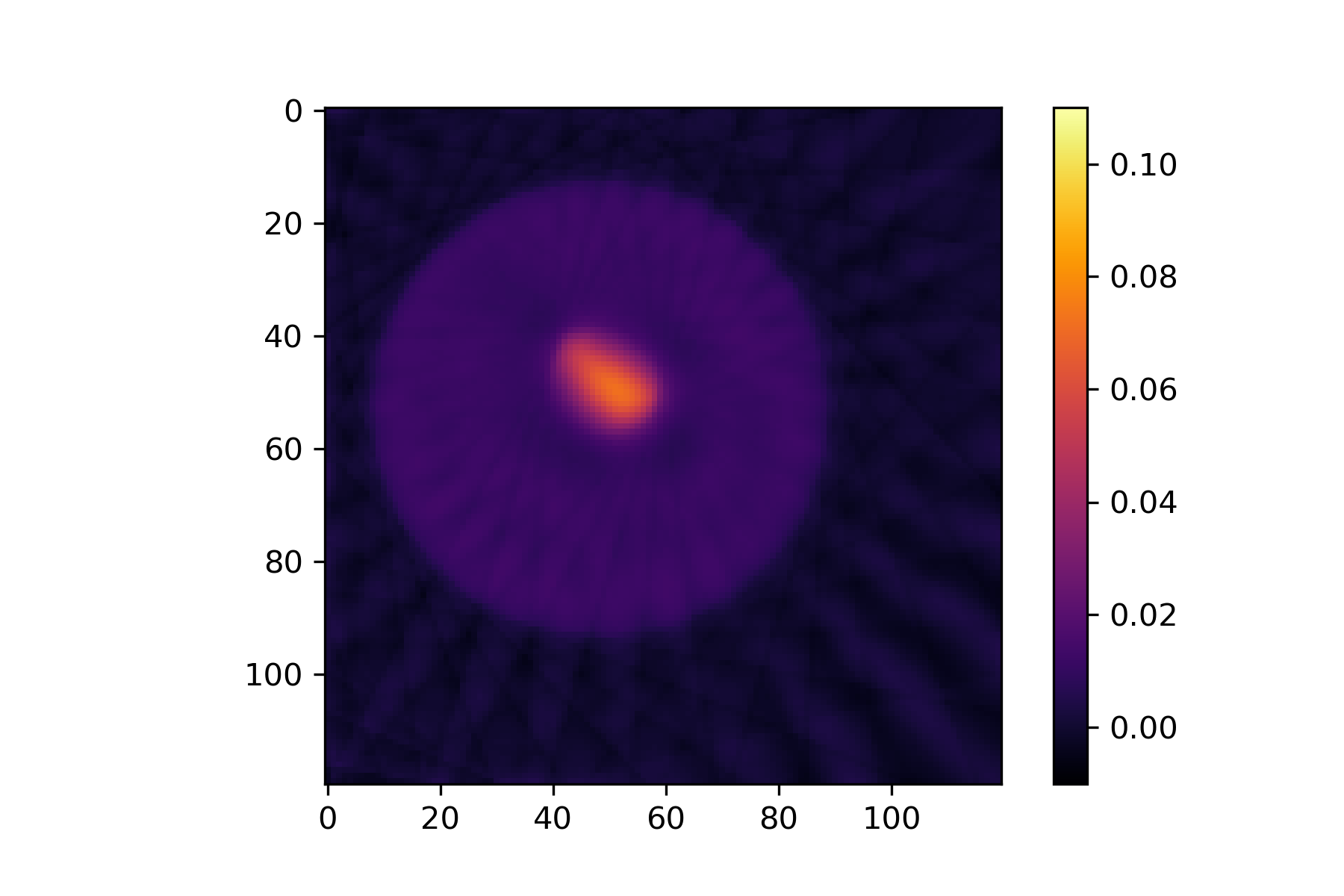}
\includegraphics[height=4.1cm,clip,trim=4.55cm 1.3cm 3.8cm 1.25cm]{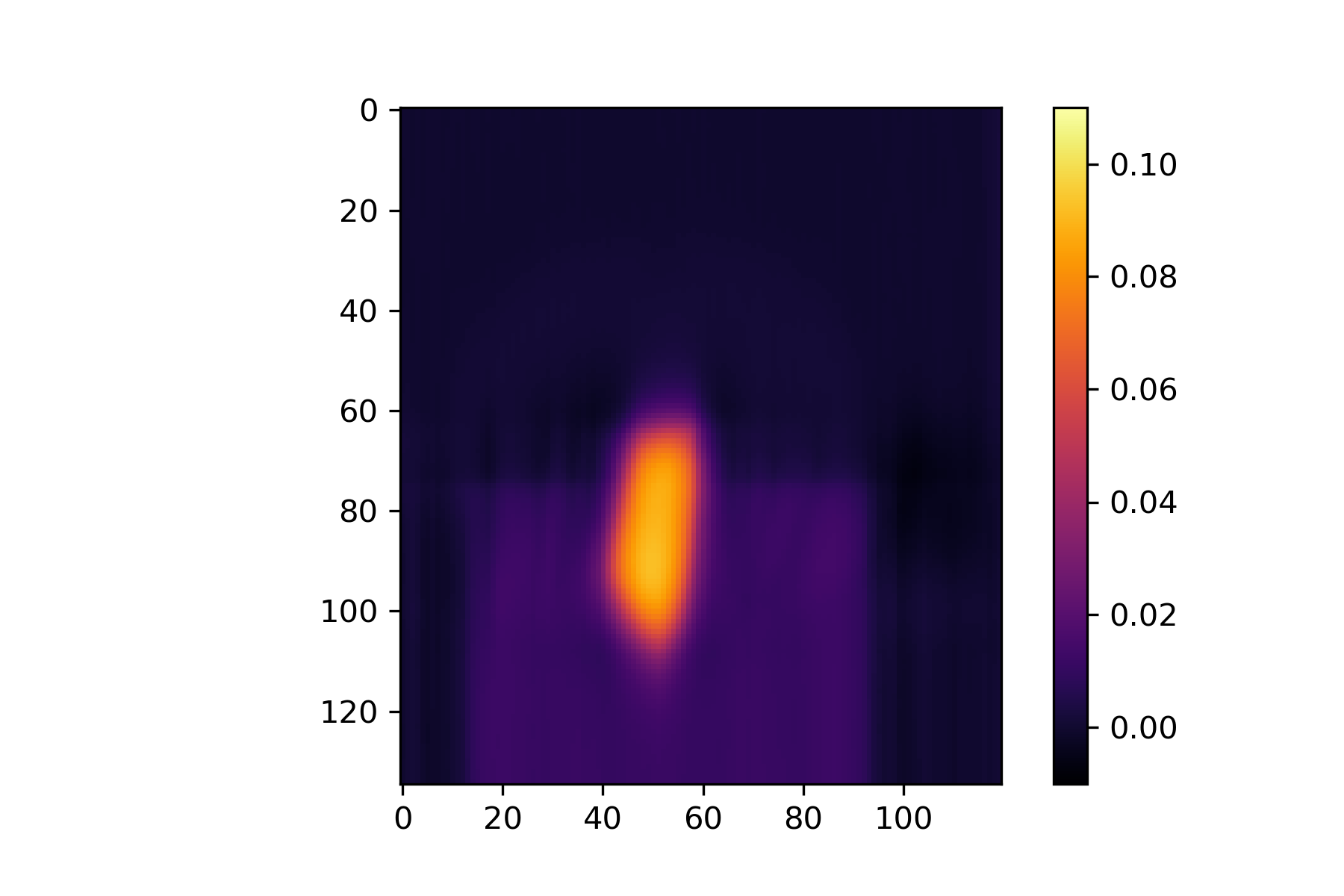}
\caption{Anisotropic Tikhonov reconstruction of 15-projection 3D steel-wire dataset. Left two: Horizontal and vertical slices, Tikhonov regularisation with horizontal smoothing ($\alpha_x = \alpha_y = 30$, $\alpha_z = 0.1$). Right two: Same, with vertical smoothing ($\alpha_x = \alpha_y = 0.1$, $\alpha_z = 60$). Colour range [-0.01,0.11].\label{fig:dlstikhonov}}
\end{figure}

In \cref{fig:dlstikhonov} (left) Tikhonov regularisation with the \code{GradientOperator} is demonstrated on the steel-wire sample. Here, $\alpha$ governs the solution smoothness similar to how the number of iterations affects CGLS solutions, with large $\alpha$ values producing smooth solutions. Here $\alpha = 1$ is used as a suitable trade-off between noise reduction and smoothing.

The block structure provides the machinery to experiment with different amounts or types of regularisation in individual dimensions in a Tikhonov setting. We consider the problem
\begin{equation}
\im^\star = \argmin_\im \left\{ \|\sysmat \im - \sino \|_2^2 + \alpha_x^2 \|D_x\im\|_2^2 + \alpha_y^2 \|D_y\im\|_2^2 + \alpha_z^2 \|D_z\im\|_2^2 \right\},
\end{equation}
where we have different regularising operators $D_x$, $D_y$, $D_z$ in each dimension and associated regularisation parameters $\alpha_x$, $\alpha_y$, $\alpha_z$. We can write this as the following block least squares problem which can be solved by CGLS:
\begin{equation}
\im^\star = \argmin_\im  \left\|
\begin{pmatrix}
\sysmat \\
\alpha_x D_x \\
\alpha_y D_y \\
\alpha_z D_z \\
\end{pmatrix}
 \im -  
 \begin{pmatrix}
\sino \\
0_x \\
0_y  \\
0_z  \\
\end{pmatrix}
 \right\|_2^2,
\end{equation}
where $0_x$, $0_y$ and $0_z$ represent zero vectors of appropriate size.

In \cref{fig:dlstikhonov} we show results for $D_x$, $D_y$ and $D_z$ being finite-difference operators in each direction, achieved by the \code{FiniteDifferenceOperator}. We show two choices of sets of regularisation parameters, namely $\alpha_x = \alpha_y = 30$, $\alpha_z = 0.1$ and $\alpha_x = \alpha_y = 0.1$, $\alpha_z = 60$. We see in the former case a large amount of smoothing occurs in the horizontal dimensions due to the larger $\alpha_x$ and $\alpha_y$ parameters, and little in the vertical dimension,  so horizontal edges are preserved. In the latter case, opposite observations can be made.

Such anisotropic regularization could be useful with objects having a layered or fibrous structure  
or if the measurement setup provides different resolution or noise properties in different dimensions, e.g., for non-standard scan trajectories such as tomosynthesis/laminography.

\subsection{Smooth convex optimisation}

\begin{table}[tb]
\centering
\caption{\fun{}s in CIL.}
\label{tab:functions}
\begin{tabular}{ll}
\hline
\textbf{Name} & \textbf{Description}\\
\hline
%
BlockFunction & Separable sum of multiple functions\\
ConstantFunction & Function taking the constant value\\
OperatorCompositionFunction & Compose function $f$ and operator $A$: $f(Ax)$\\
IndicatorBox & Indicator function for box (lower/upper) constraints \\
KullbackLeibler & Kullback-Leibler divergence data fidelity\\
L1Norm & $L^1$-norm: $\|x\|_1 = \sum_i |x_i|$ \\
L2NormSquared & Squared $L^2$-norm: $\|x\|_2^2 = \sum_i x_i^2$ \\
LeastSquares & Least-squares data fidelity: $\|Ax-b\|_2^2$\\
MixedL21Norm & Mixed $L^{2,1}$-norm: $\| (U_1; U_2)\|_{2,1} = \| (U_1^2 + U_2^2)^{1/2} \|_1$ \\
SmoothMixedL21Norm & Smooth $L^{2,1}$-norm: $\| (U_1; U_2)\|_{2,1}^\text{S} = \| (U_1^2 + U_2^2 + \beta^2)^{1/2} \|_1$\\
WeightedL2NormSquared & Weighted squared $L^2$-norm: $\|x\|_w^2 = \sum_i (w_i \cdot x_i^2)$  \\ \hline
\end{tabular}
\vspace*{-4pt}
\end{table}

CIL supports the formulation and solution of more general  optimisation problems. One problem class supported is unconstrained smooth convex optimisation problems, 
\begin{equation}
\im^\star = \argmin_\im f(\im).
\end{equation}
Here $f$ is a differentiable, convex, so-called $L$-smooth function, that is its gradient $\nabla f$ is $L$-Lipschitz continuous: $\|\nabla f(\im_1) - \nabla f(\im_2) \|_2 \leq L \| \im_1 - \im_2\|_2,\; \forall \im_1, \im_2$ for some $L>0$ referred to as the Lipschitz parameter. CIL represents functions by the \fun{} class, which maps an \code{ImageData} or \code{AcquisitionData} to a real number. Differentiable functions provide the method \code{gradient} to allow first-order optimisation methods to work; at present CIL provides a Gradient Descent method \code{GD} with a constant or back-tracking line search for step size selection. CIL \fun supports algebra so the user can formulate for example linear combinations of \fun objects and solve with the \code{GD} algorithm.

As example we can formulate and solve the Tikhonov problem \cref{eq:tikhonov} with \code{GD} as
\begin{center}
\begin{tcolorbox}[
    enhanced,
    attach boxed title to top center={yshift=-2mm},
    colback=darkspringgreen!20,
    colframe=darkspringgreen,
    colbacktitle=darkspringgreen,
    title=Set up and run Gradient Descent for Tikhonov regularisation,
    text width = 15cm,
    fonttitle=\bfseries\color{white},
    boxed title style={size=small,colframe=darkspringgreen,sharp corners},
    sharp corners,
]
\begin{minted}{python}
f1 = LeastSquares(A, b)
f2 = OperatorCompositionFunction(L2NormSquared(), L)
f  = f1 + (alpha**2)*f2
myGD = GD(initial=x0, objective_function=f)
myGD.run(1000, verbose=1)
\end{minted}
\end{tcolorbox}
\end{center}
%
%
Here \code{LeastSquares(A,b)}, representing $\|A \cdot -b\|_2^2$, and \code{L2NormSquared}, representing $\|\cdot\|_2^2$, are examples from the \fun{} class. 
With \code{OperatorCompositionFunction} a function can be composed with an operator, here \code{L}, to form a composite function $\|L\cdot\|_2^2$.
An overview of \fun{} types currently in CIL is provided in \cref{tab:functions}.
Another example using a smooth approximation of non-smooth total variation regularisation will be given in \cref{sec:neutron}.

\subsection{Non-smooth convex optimisation with simple proximal mapping}

\begin{figure}[t]
\includegraphics[height=4.1cm,clip,trim=3.7cm 1.3cm 3.9cm 1.25cm]{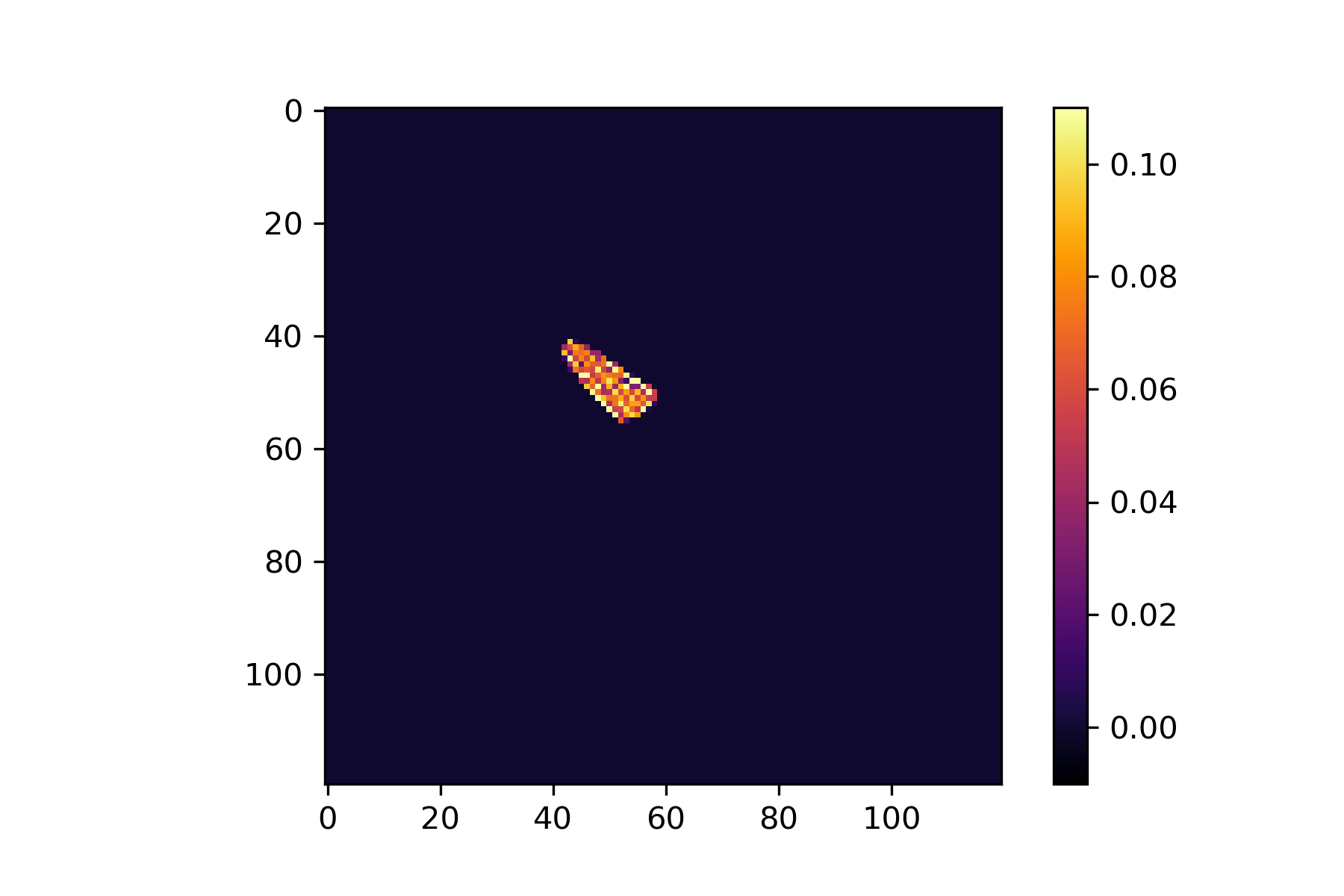}
\includegraphics[height=4.1cm,clip,trim=4.55cm 1.3cm 3.8cm 1.25cm]{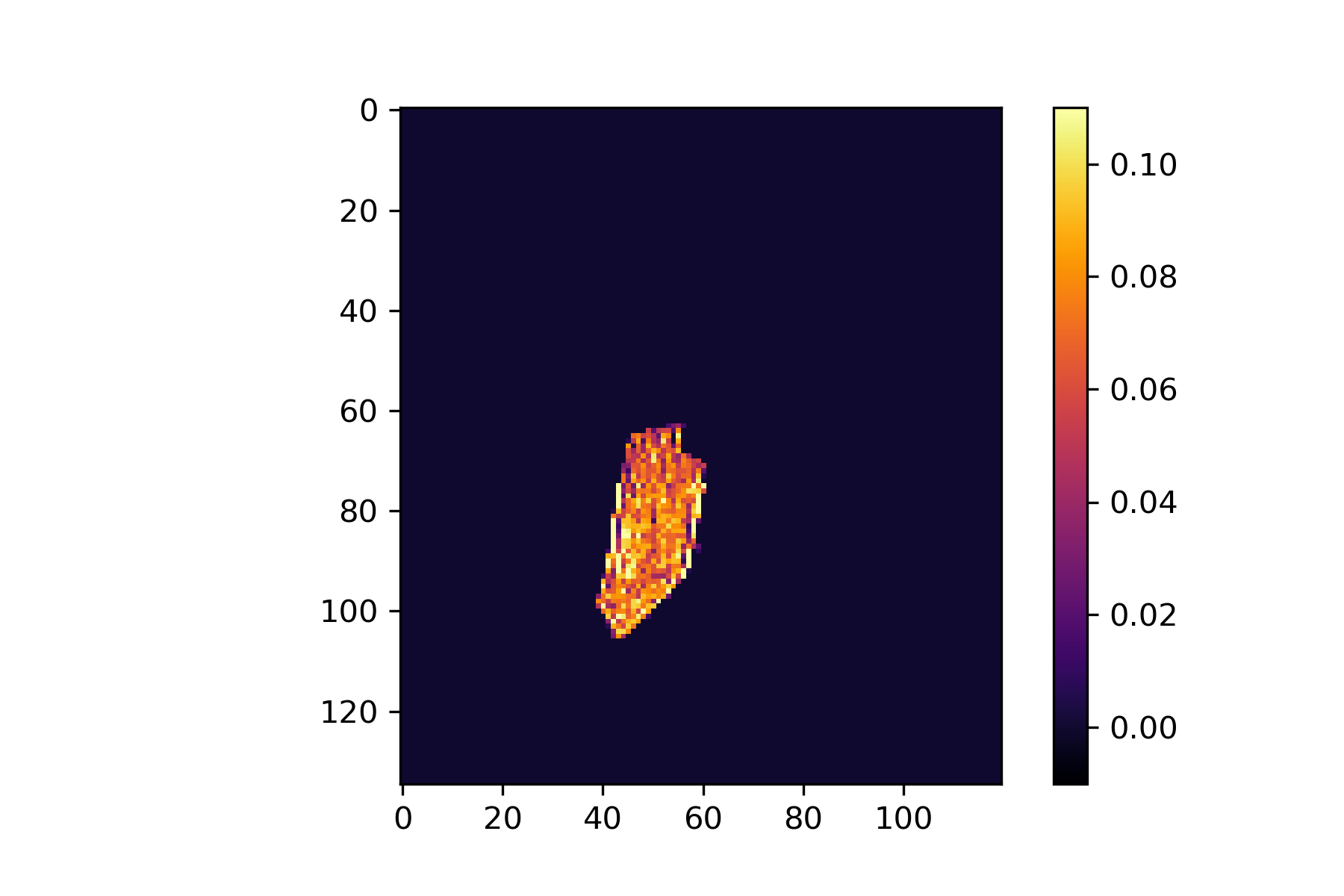}
\includegraphics[height=4.1cm,clip,trim=3.7cm 1.3cm 3.9cm 1.25cm]{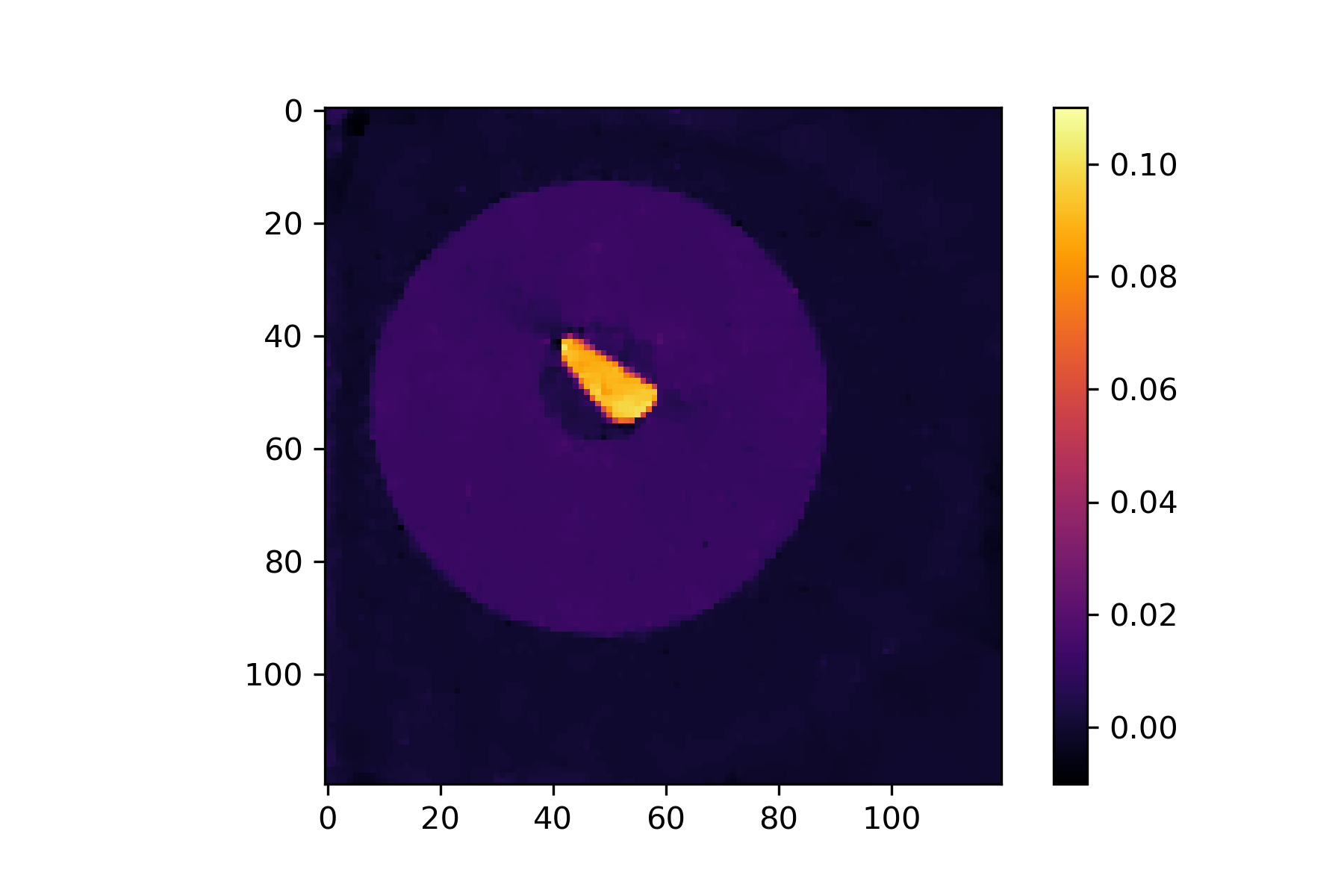}
\includegraphics[height=4.1cm,clip,trim=4.55cm 1.3cm 3.8cm 1.25cm]{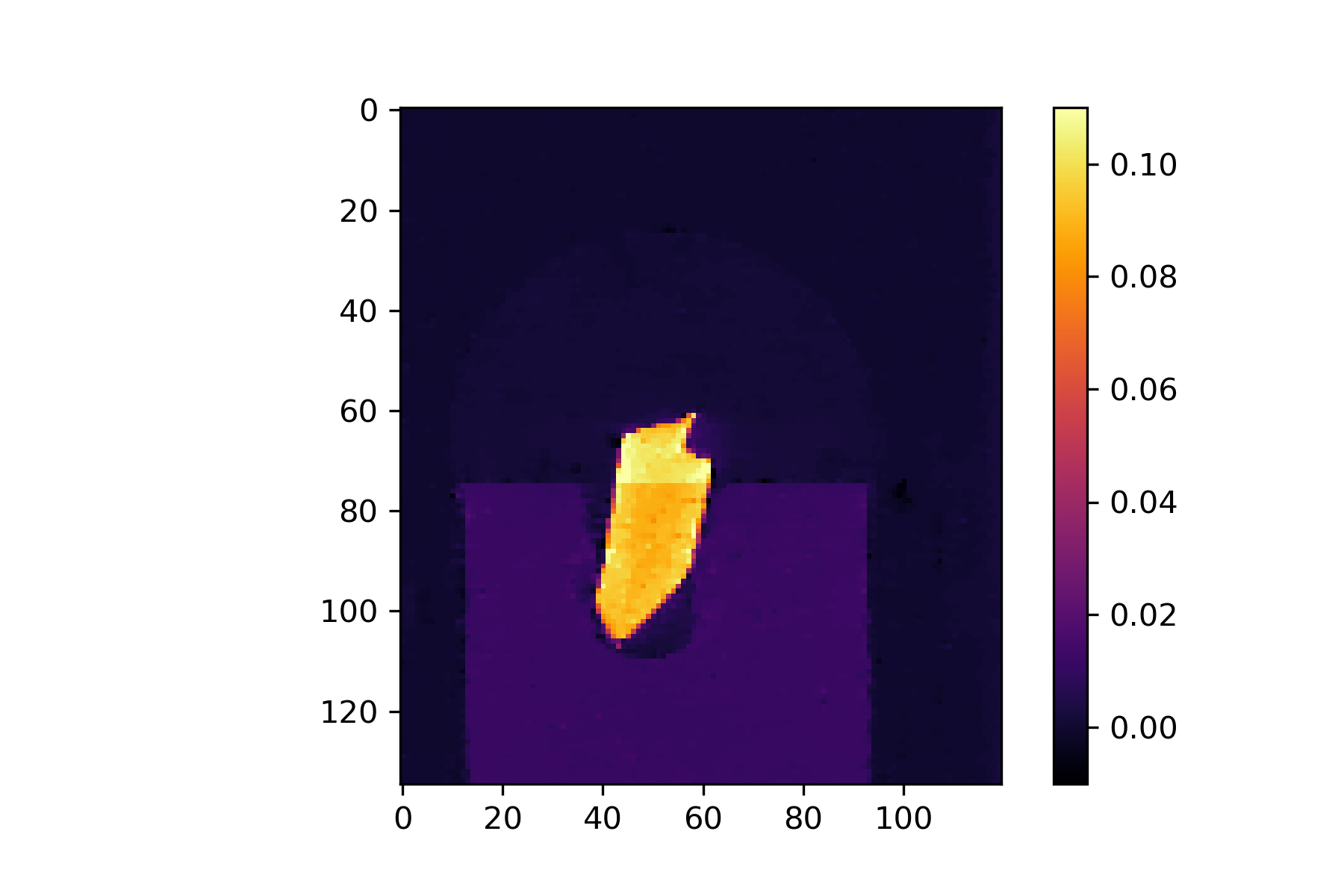}
\caption{FISTA reconstruction of 15-projection 3D steel-wire dataset. Left two: $L^1$-norm regularization with large regularisation parameter of $\alpha=30$ forces all pixels but in steel wire to zero. Right two: TV-regularisation with $\alpha=0.02$ removes streaks and noise and preserves edges. Colour range [-0.01,0.11].\label{fig:dlsfista}}
\end{figure}

Many useful reconstruction methods are formulated as non-smooth optimisation problems. Of specific interest in recent years has been sparsity-exploiting regularisation such as the $L^1$-norm and total variation (TV). TV-regularisation for example has been shown capable of producing high-quality images from severely undersampled data whereas FBP produces highly noisy, streaky images. A particular problem class of interest can be formulated as
\begin{equation}
\im^\star = \argmin_{\im} \Big\{ f(\im) + g(\im) \Big\},  \label{eq:fistaproblem}
\end{equation}
where $f$ is $L$-smooth and $g$ may be non-smooth. This problem can be solved by the Fast Iterative Shrinkage-Thresholding Algorithm (FISTA) \cite{FISTASIIMS,FISTAIEEE}, which is available in CIL as \code{FISTA}. FISTA makes use of $f$ being smooth by calling \code{f.gradient} and assumes for $g$ that the so-called proximal mapping,
\begin{align}
\mathrm{prox}_{\tau g}(u) = \argmin_v \left\{\tau g(v) + \frac{1}{2}\|v-u\|_2^2 \right\}
\end{align}
for a positive parameter $\tau$ is available as \code{g.proximal}. This means that FISTA is useful when $g$ is ``proximable'', i.e., where an analytical expression for the proximal mapping exists, or it can be computed efficiently numerically. 

A simple, but useful case, for FISTA is to enforce constraints on the solution, i.e.,  require $\im \in C$, where $C$ is a convex set. In this case $g$ is set to the (convex analysis) indicator function of $C$, i.e.,
\begin{equation}
\iota_C(\im) =  
\left\{\begin{array}{lll}
0 & \text{if} & \im \in C \\
\infty & \text{else}. &
\end{array}  \right.
\end{equation}
The proximal mapping of an indicator function is simply a projection onto the convex set; for simple lower and upper bound constraints this is provided in CIL as \code{IndicatorBox}. FISTA with non-negativity constraints is achieved with the following lines of code:
\begin{center}
\begin{tcolorbox}[
    enhanced,
    attach boxed title to top center={yshift=-2mm},
    colback=darkspringgreen!20,
    colframe=darkspringgreen,
    colbacktitle=darkspringgreen,
    title=Set up and run FISTA for non-negative least-squares problem,
    text width = 15cm,
    fonttitle=\bfseries\color{white},
    boxed title style={size=small,colframe=darkspringgreen,sharp corners},
    sharp corners,
]
\begin{minted}{python}
F = LeastSquares(A, b)
G = IndicatorBox(lower=0.0) 
myFISTA = FISTA(f=F, g=G, initial=x0, max_iteration=1000)
myFISTA.run(300, verbose=1)
\end{minted}
\end{tcolorbox}
\end{center}
Another simple non-smooth case is $L^1$-norm regularisation, i.e., using $\|\im\|_1 = \sum_j |\im_j|$ as regulariser. This is non-differentiable at $0$ and a closed-form expression for the proximal mapping is known as the so-called soft-thresholding. In CIL this is available as \code{L1Norm} and can be achieved with the same code, only with the second line replaced by
\begin{center}
\begin{tcolorbox}[
    enhanced,
    attach boxed title to top center={yshift=-2mm},
    colback=darkspringgreen!20,
    colframe=darkspringgreen,
    colbacktitle=darkspringgreen,
    title=Set up L1 regulariser for use in FISTA,
    text width = 15cm,
    fonttitle=\bfseries\color{white},
    boxed title style={size=small,colframe=darkspringgreen,sharp corners},
    sharp corners,
]
\begin{minted}{python}
alpha = 100
G = alpha*L1Norm()
\end{minted}
\end{tcolorbox}
\end{center}
The resulting steel-wire dataset reconstruction is shown in \cref{fig:dlsfista}.

FISTA can also be used whenever a numerical method is available for the proximal mapping of $g$; one such case is the (discrete, isotropic) Total Variation (TV). TV is the mixed $L^{2,1}$-norm of the gradient image, 
\begin{equation}
g_\text{TV}(u) = \|D u \|_{2,1} = \left\|\binom{D_x}{D_y} u\right\|_{2,1} = \left\|\sqrt{ (D_x u)^2 + (D_y u)^2}\right\|_1 , \label{eq:TVdef}
\end{equation}
where $D = (D_x; D_y)$ is the gradient operator as before and 
the $L^2$-norm combines the $x$ and $y$ differences before the $L^1$-norm sums over all voxels.
CIL implements this in \code{TotalVariation} using the FGP method from \cite{FISTAIEEE}. Using the FISTA code above we can achieve this with 
\begin{center}
\begin{tcolorbox}[
    enhanced,
    attach boxed title to top center={yshift=-2mm},
    colback=darkspringgreen!20,
    colframe=darkspringgreen,
    colbacktitle=darkspringgreen,
    title=Set up TV regulariser for use in FISTA,
    text width = 15cm,
    fonttitle=\bfseries\color{white},
    boxed title style={size=small,colframe=darkspringgreen,sharp corners},
    sharp corners,
]
\begin{minted}{python}
alpha = 0.02
G = alpha*TotalVariation()
\end{minted}
\end{tcolorbox}
\end{center}
The resulting reconstruction is shown in \cref{fig:dlsfista} and clearly demonstrates the edge-preserving, noise-reducing and streak-removing capabilities of TV-regularisation.

\subsection{Non-smooth convex optimisation using splitting methods}

\begin{figure}[t]
\includegraphics[height=4.5cm,clip,trim=3.7cm 1.3cm 3.9cm 1.25cm]{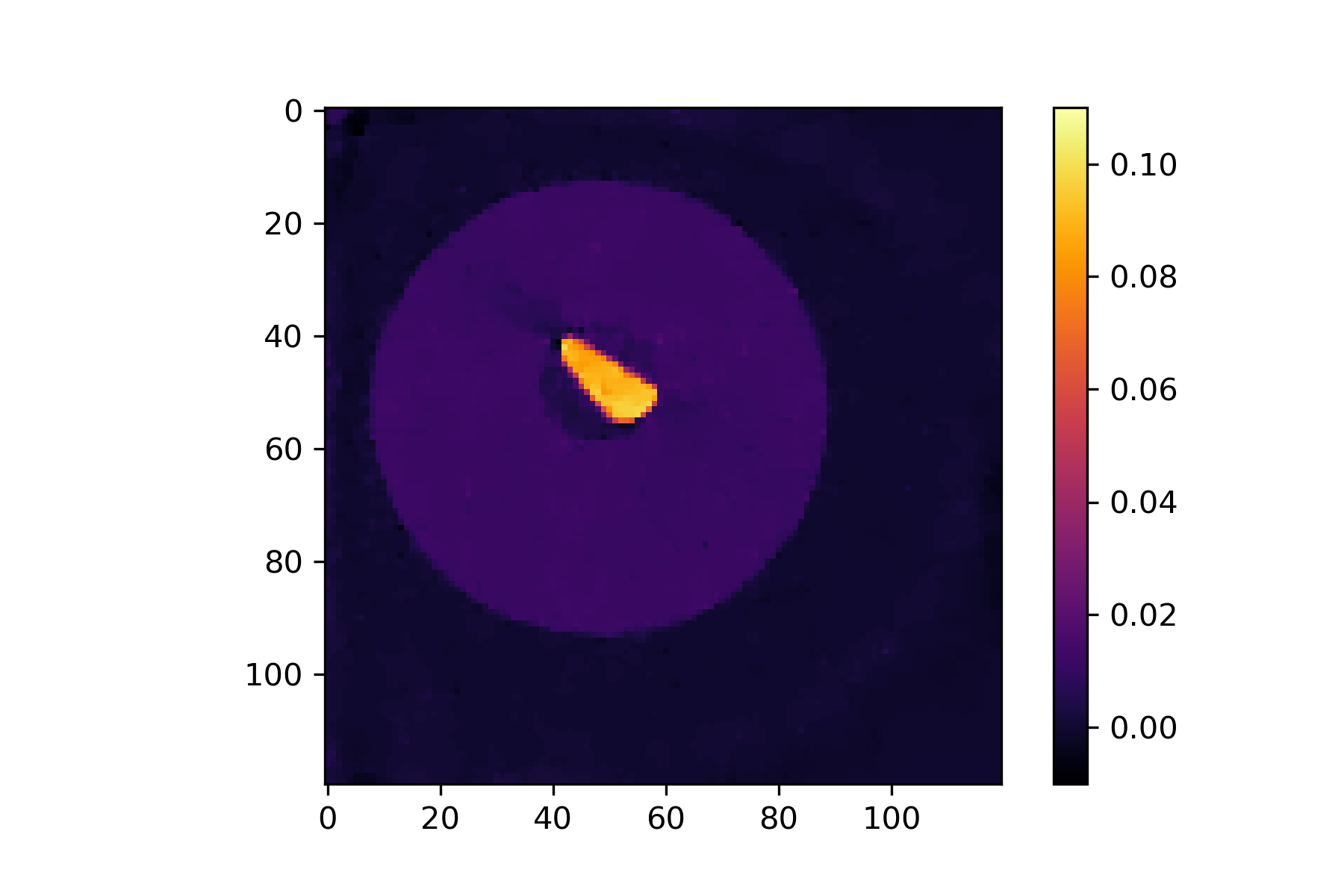}
\includegraphics[height=4.5cm,clip,trim=4.55cm 1.3cm 3.8cm 1.25cm]{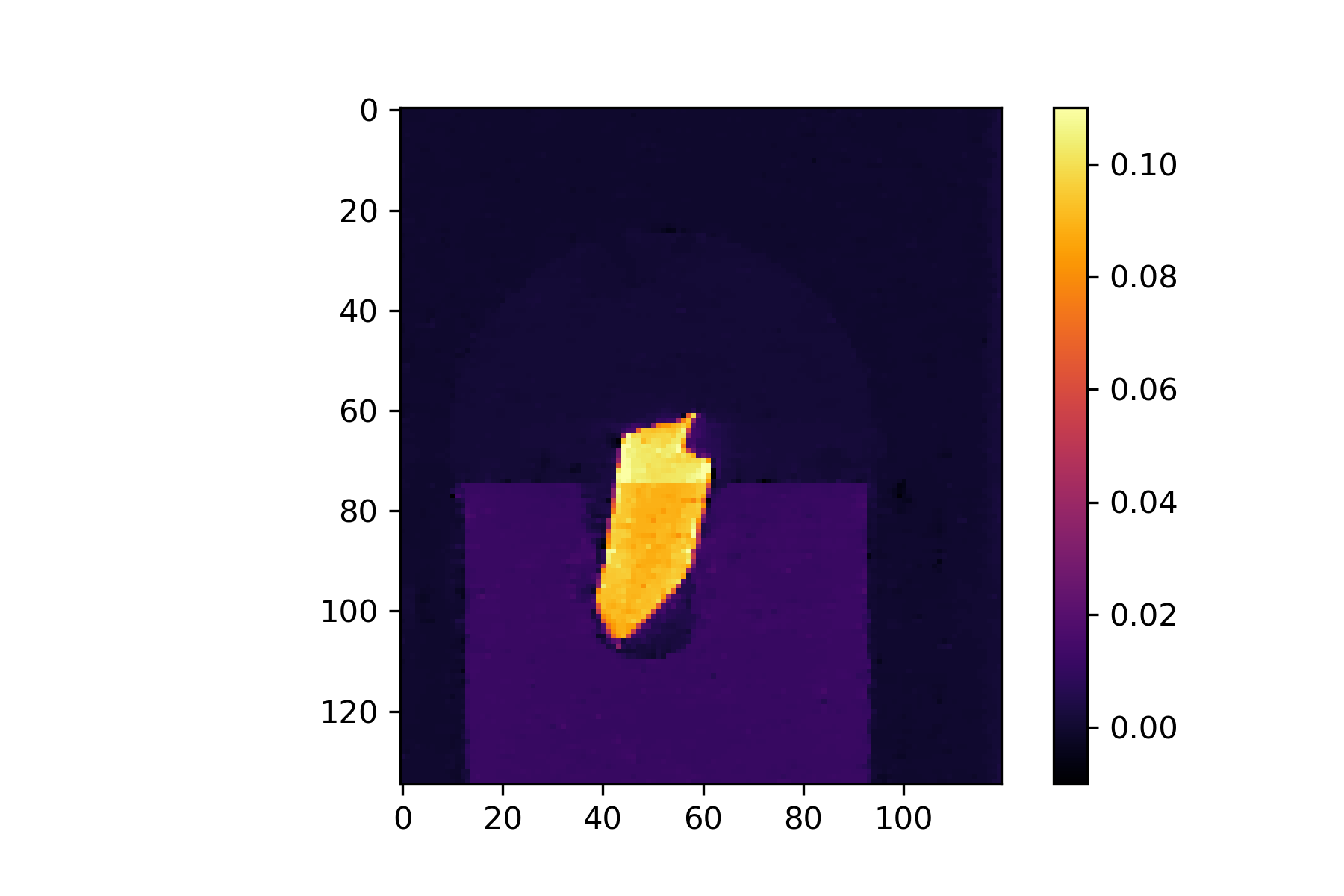}
\;
\includegraphics[width=0.43\linewidth,clip,trim=0.4cm 0.2cm 1.5cm 0.75cm]{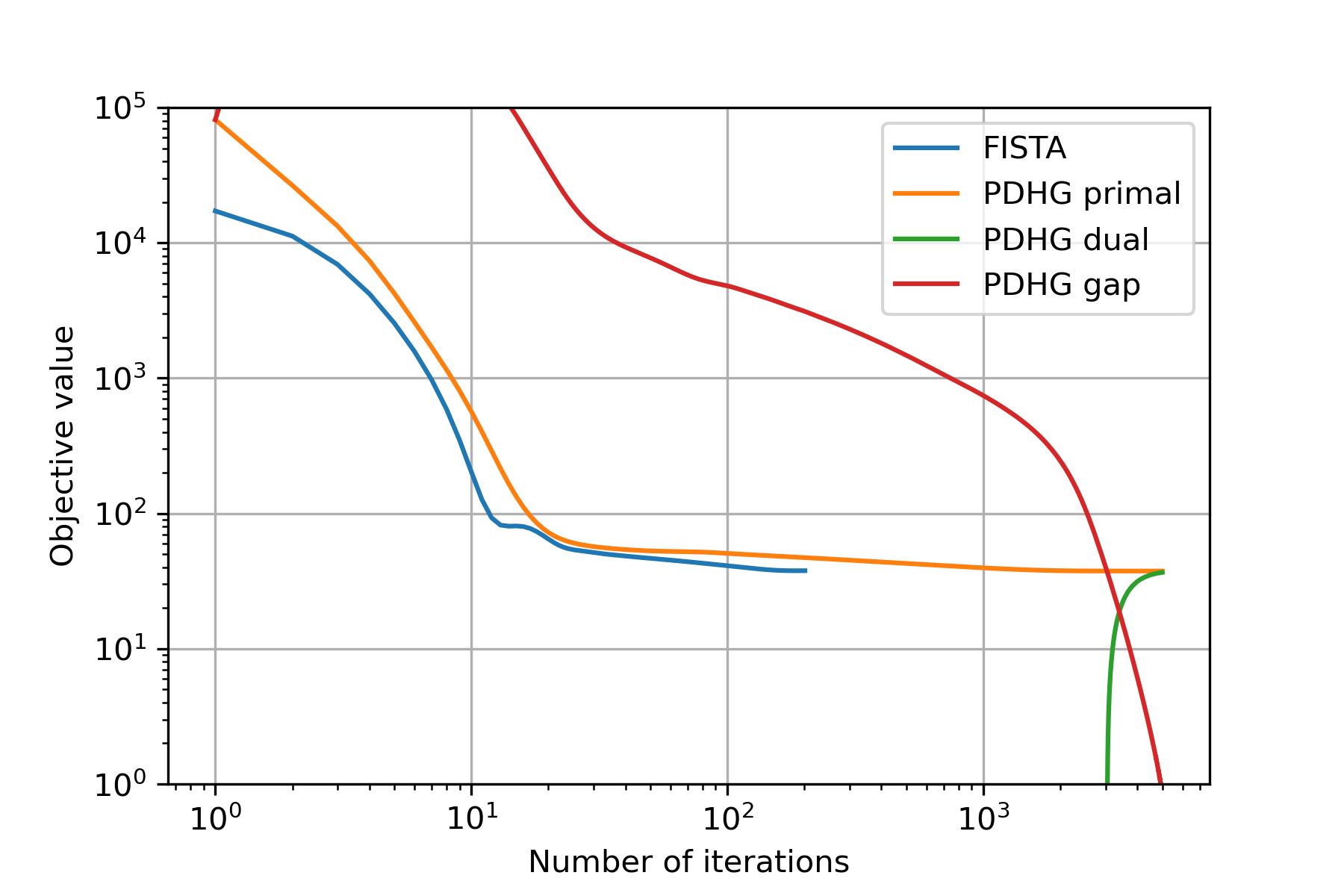}
\caption{PDHG reconstruction of 15-projection 3D steel-wire dataset. Left two:  TV-regularisation with $\alpha=0.02$, reproduces same result as FISTA in \cref{fig:dlsfista} on same case and parameter choice, thus validating algorithms against each other. Colour range [-0.01,0.11]. Right: Objective value histories (log-log) for FISTA and PDHG on TV-regularisation problem. Both algorithms reach the same (primal) objective value, FISTA taking fewer but slower iterations. The primal-dual gap for PDHG (difference between primal and dual objectives) approaches zero indicating convergence.
\label{fig:tvpdhg}}
\end{figure}

When the non-smooth function is not proximable, we may consider so-called \emph{splitting methods} for solving a more general class of problems, namely 
\begin{equation}
\im^\star = \argmin_\im \Big\{ f(K\im) + g(\im) \Big\}, \label{eq:PDHGproblem}
\end{equation}
where $f$ and $g$ are convex (possibly) non-smooth functions and $K$ a linear operator. The key change from the FISTA problem is the \emph{splitting} of the complicated $f(K(u))$, which as a whole may \emph{not} be proximable, into simpler parts $f$ and $K$ to be handled separately. CIL provides two algorithms for solving this problem, depending on properties of $f$ and assuming that $g$ is proximable. If $f$ is proximable, then the linearised ADMM method  \cite{ParikhBoyd2014} can be used; available as \code{LADMM} in CIL. If the so-called convex conjugate, $f^*$, of $f$ is proximable, then the Primal Dual Hybrid Gradient (PDHG) method \cite{Esser2010, ChambollePock2011, Sidky_2012}, also known as the Chambolle-Pock method, may be used; this is known as \code{PDHG} in CIL.

\begin{table}[tb]
\centering
\caption{\alg{}s in CIL.}
\label{tab:algorithms}
\begin{tabular}{lll}
\hline
\textbf{Name} & \textbf{Description} & \textbf{Problem type solved}  \\
\hline
CGLS & Conjugate Gradient Least Squares & Least squares\\
SIRT & Simultaneous Iterative Reconstruction Technique & Weighted least squares \\
GD & Gradient Descent & Smooth\\
FISTA & Fast Iterative Shrinkage-Thresholding Algorithm & Smooth + non-smooth\\
LADMM & Linearised Alternating Direction Method of Multipliers & Non-smooth\\
PDHG & Primal Dual Hybrid Gradient & Non-smooth\\
SPDHG & Stochastic Primal Dual Hybrid Gradient & Non-smooth \\\hline
\end{tabular}
\vspace*{-4pt}
\end{table}

In fact an even wider class of problems can be handled using this formulation, namely
\begin{equation}
\im^\star = \argmin_\im \left\{ \sum_i f_i(K_i\im) + g(\im) \right\},
\end{equation}
i.e., where the composite function $f(K\cdot)$ can be written as a separable sum
\begin{equation}
 f(K\im) = \sum_i f_i(K_i\im).
\end{equation}
In CIL we can express such a function using a \code{BlockOperator}, as also used in the Tikhonov example, and a \code{BlockFunction}, which essentially holds a list of \fun{} objects.

Here we demonstrate this setup by using \code{PDHG} to solve the TV-regularised least-squares problem. 
As shown in \cite{Sidky_2012} this problem can be written in the required form as
\begin{equation}
f = \binom{f_1}{f_2} = \binom{ \|\cdot -b \|_2^2}{\alpha \|\cdot\|_{2,1}}, \qquad
K = \binom{A}{D},\qquad 
g(u) = 0.
\end{equation}
In CIL this can be written succinctly as (with a specific choice of regularisation parameter):

%
%
%
\begin{center}
\begin{tcolorbox}[
    enhanced,
    attach boxed title to top center={yshift=-2mm},
    colback=darkspringgreen!20,
    colframe=darkspringgreen,
    colbacktitle=darkspringgreen,
    title=Set up and run PDHG for TV-regularised least-squares problem,
    text width = 15cm,
    fonttitle=\bfseries\color{white},
    boxed title style={size=small,colframe=darkspringgreen,sharp corners},
    sharp corners,
]
\begin{minted}{python}
alpha = 0.02
F = BlockFunction(L2NormSquared(b=b), alpha*MixedL21Norm())
K = BlockOperator(A, GradientOperator(ig))
G = ZeroFunction()
myPDHG = PDHG(f=F, operator=K, g=G, max_iteration=10000)
myPDHG.run(5000, verbose=2)
\end{minted}
\end{tcolorbox}
\end{center}

\Cref{fig:tvpdhg} shows the resulting steel-wire dataset reconstruction which appears identical to the result of FISTA on the same problem (\cref{fig:dlsfista}), and as such validates the two algorithms against each other.

CIL \alg{}s have the option to save the history of objective values so the progress and convergence can be monitored. PDHG is a primal-dual algorithm, which means that the so-called dual maximisation problem of \cref{eq:PDHGproblem}, which is referred to as the primal problem, is solved simultaneously. In \code{PDHG} the dual objective values are also available. The primal-dual gap, which is the difference between the primal and dual objective values, is useful for monitoring convergence as it should approach zero when the iterates converge to the solution.

\Cref{fig:tvpdhg} (right) compares the primal objective, dual objective and primal-dual gap history with the objective history for FISTA on the same problem. The (primal) objectives settle at roughly the same level, again confirming that the two algorithms achieve essentially the same solution. FISTA used fewer iterations, but each iteration took about 25 times as long as a PDHG iteration. The dual objective is negative until around 3000 iterations, and the primal-dual gap is seen to approach zero, thus confirming convergence. CIL makes such algorithm comparisons straightforward. It should be stressed that the particular convergence behavior observed for \code{FISTA} and \code{PDHG} depends on internal algorithm parameters such as step sizes for which default values were used here. The user may experiment with tuning these parameters to obtain faster convergence, for example for \code{PDHG} the primal and dual step sizes may be set using the inputs \code{sigma} and \code{tau}.

In addition to PDHG a stochastic variant SPDHG \cite{SPDHG} that can sometimes accelerate reconstruction substantially by working on problem subsets is provided in CIL as \code{SPDHG}; this is demonstrated in the Part II article \cite{CIL2}.

An overview of all the algorithms currently supplied by CIL is provided in \cref{tab:algorithms}.


\section{Exemplar studies using CIL} \label{sec:examples}
 
This section presents 3 illustrative examples each demonstrating different functionality of CIL. All code and data to reproduce the results are provided, see Data Accessibility.

\subsection{Neutron tomography with golden-angle data} \label{sec:neutron}

This example demonstrates how CIL can handle other imaging modalities than X-ray, a  non-standard scan geometry, and easily compare reconstruction algorithms.

Contrary to X-rays, neutrons interact with atomic nuclei rather than electrons that surround them, which yields a different contrast mechanism, 
e.g., for neutrons hydrogen is highly attenuating while lead is almost transparent.
Nevertheless, 
neutron data can be modelled with the Radon transform and reconstructed with the same techniques as X-ray data.

\begin{figure}[tb]
\includegraphics[width=0.335\linewidth]{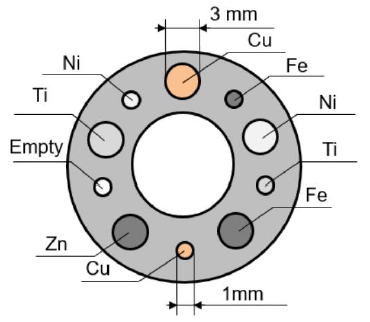}\,\,
\includegraphics[width=0.335\linewidth]{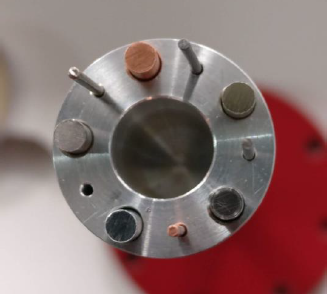}\,\,
\includegraphics[width=0.2975\linewidth]{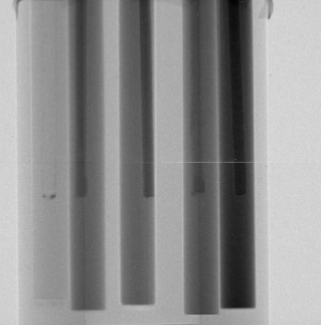}\\
\includegraphics[trim=0.5cm 1.3cm 1.5cm 1.5cm,clip,width=0.535\linewidth]{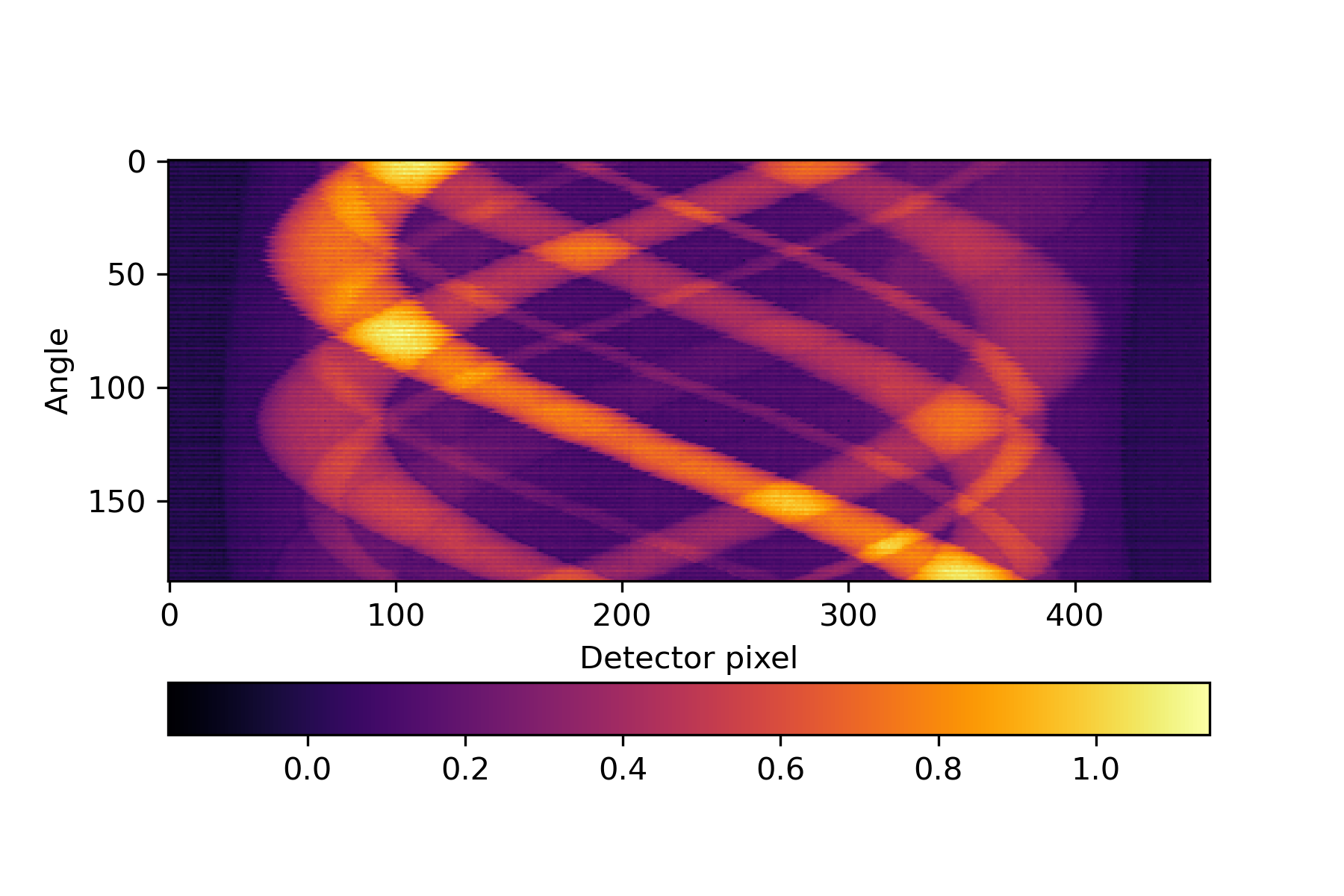}\;\,
\includegraphics[trim=0.2cm 0.2cm 1.0cm 1cm,clip,width=0.46\linewidth]{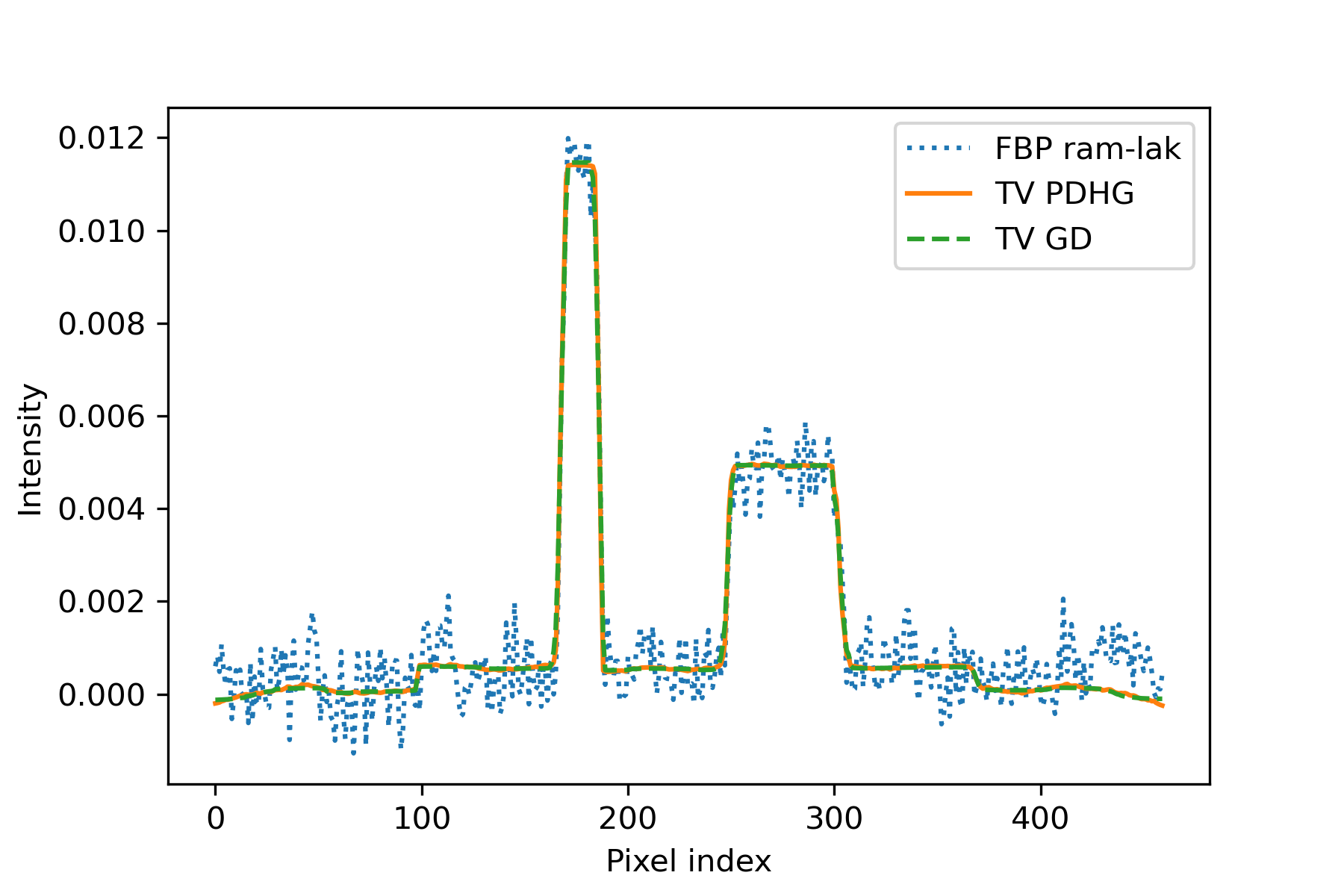}\\[-0.3cm]

\includegraphics[trim=3.7cm 1.3cm 3.8cm 1.2cm,clip,width=0.33\linewidth]{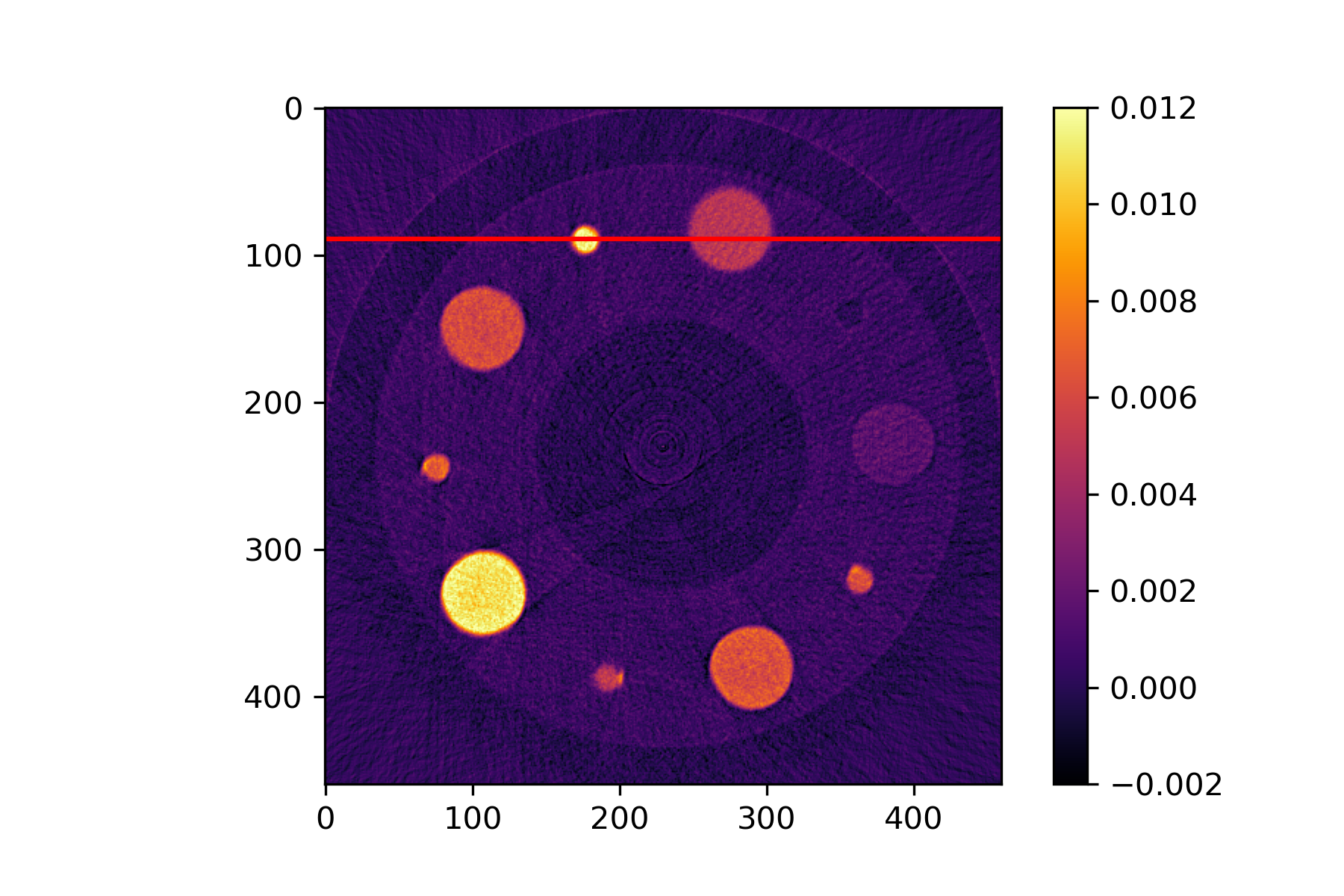}
\includegraphics[trim=3.7cm 1.3cm 3.8cm 1.2cm,clip,width=0.33\linewidth]{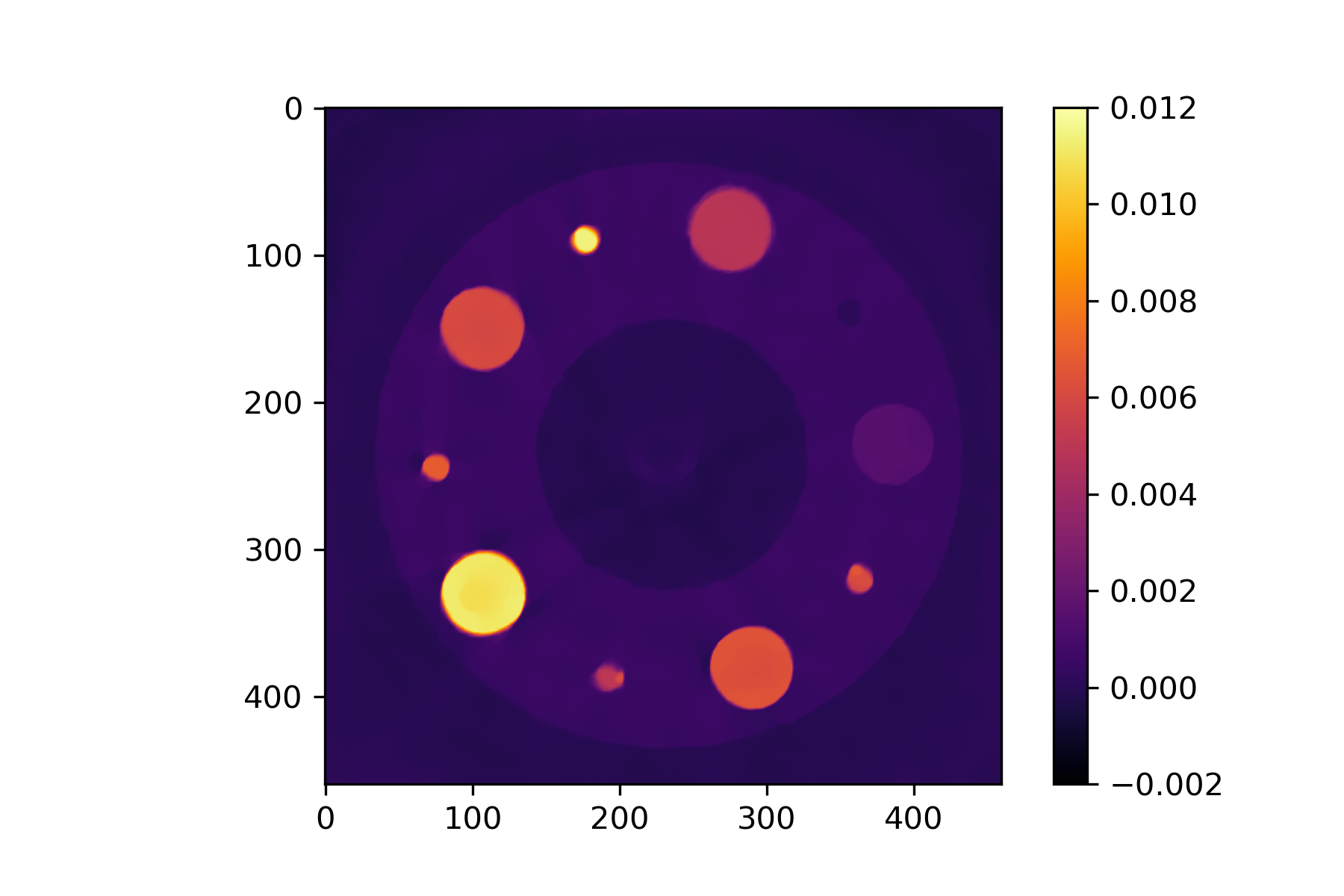}
\includegraphics[trim=3.7cm 1.3cm 3.8cm 1.2cm,clip,width=0.33\linewidth]{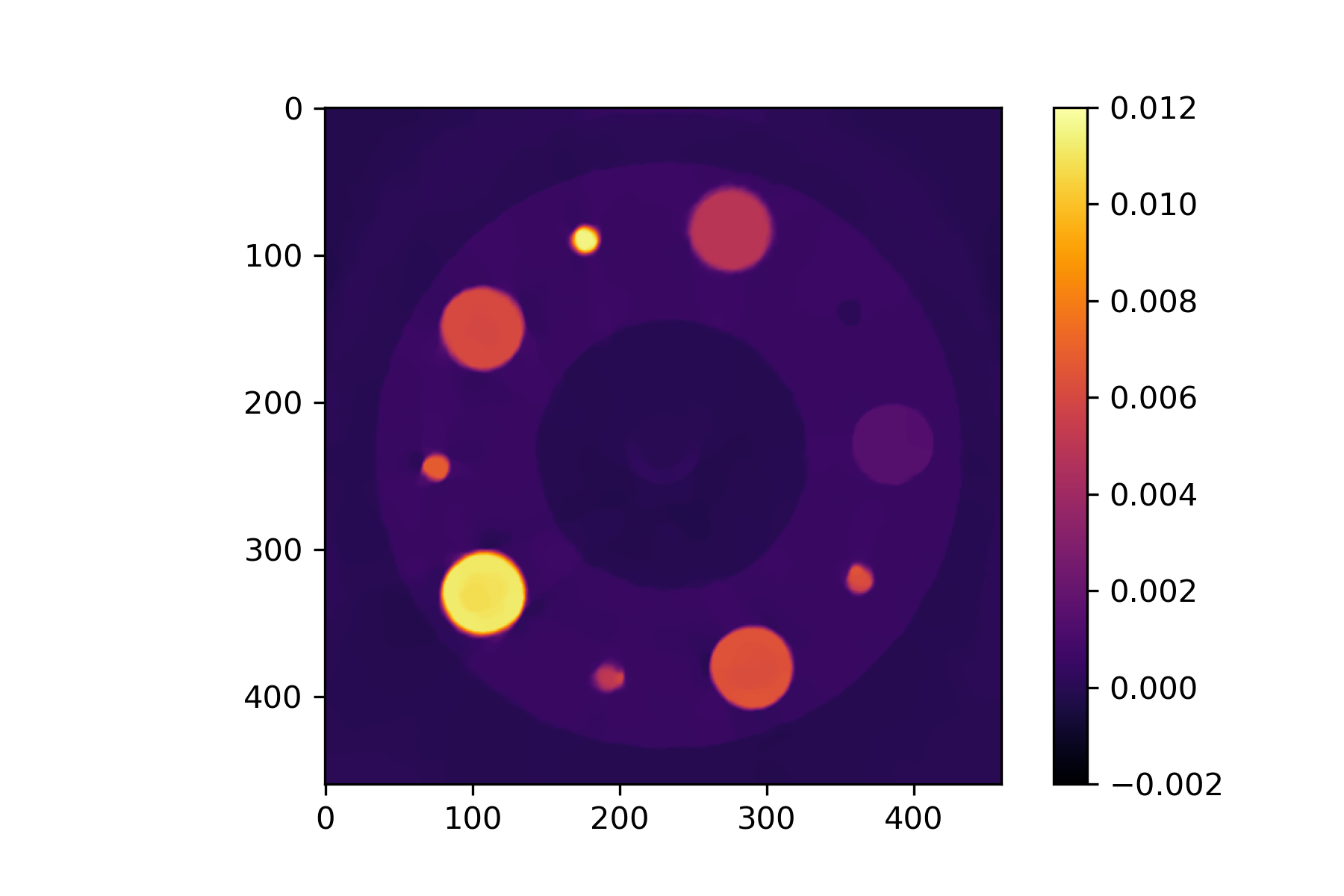}
\caption{IMAT neutron tomography dataset. Top row: (left) top-view schematic of high-purity elemental metal rod sample; (centre) top-view photograph; (right) single raw projection image showing rods of different absorption. Middle row: (left) preprocessed slice sinogram; (right) horizontal line profile of FBP, PDHG TV and GD TV reconstruction along line shown on image below. Bottom row: (left) slice reconstructions, FBP; (centre) TV reconstruction with \code{PDHG};  (right) STV reconstruction with \code{GD}. Colour range [-0.002, 0.012]. }
\label{fig:imat}
\end{figure}

A benchmarking neutron tomography dataset (\cref{fig:imat}) was acquired at  the IMAT beamline \cite{Burca_2013,KOCKELMANN201571} of the ISIS Neutron and Muon Source, Harwell, UK. The raw data is available at \cite{jorgensenIMATdata2019} and a processed subset for this paper is available from \cite{jorgensenIMATzenodo}. The test phantom consisted of an Al cylinder of diameter 22 mm with cylindrical holes holding 1mm and 3mm rods of high-purity elemental Cu, Fe, Ni, Ti, and Zn rods. 186 projections each 512-by-512 pixels in size 0.055 mm were acquired using the non-standard golden-angle mode \cite{Kohler2004APA} (angular steps of $\frac{1}{2} (\sqrt{5}-1) \cdot  180^\circ = 111.24...^\circ$) rather than sequential small angular increments. This was to provide complete angular coverage in case of early experiment termination and to allow experimenting with reconstruction from a reduced number of projections. An energy-sensitive micro-channel plate (MCP) detector was used \cite{mcptremsin,Tremsin_2014} providing raw data in 2332 energy bins per pixel, which were processed and summed to simulate a conventional white-beam absorption-contrast data set for the present paper. Reconstruction and analysis of a similar energy-resolved data set is given in \cite{Evelina}.

We use \code{TIFFStackReader} to load the data, several \proc{} instances to preprocess it, and initially \code{FBP} to reconstruct it. We compare with TV-regularisation, \cref{eq:TVdef}, solved with \code{MixedL21Norm} and \code{PDHG} using $\alpha=1$ and 30000 iterations, and further with a smoothed variant of TV (STV) using \code{SmoothMixedL21Norm}. The latter makes the optimisation problem smooth, so it can be solved using \code{GD}, using the same $\alpha$ and 10000 iterations.

The sinogram for a single slice is shown in \cref{fig:imat} along with FBP, TV and STV reconstructions and a horizontal line profile plot as marked by the red line. The FBP reconstruction recovers the main sample features, however it is contaminated by noise, ring artifacts and streak artifacts emanating from the highest-attenuating rods. The TV and STV reconstructions remove these artifacts, while preserving edges. We see that the STV approximates the non-smooth TV very well; this also serves to validate the reconstruction algorithms against one another.

\subsection{Non-standard acquisition: X-ray laminography} \label{sec:lami}

This example demonstrates how even more general acquisition geometries can be processed using CIL, and how \module{cil.plugins.ccpi\_regularisation} allows CIL to use GPU-accelerated implementations of regularising functions available in the CCPi-RGL toolkit \cite{Kazantsev2019RGLTK}. Furthermore, unlike the examples up to now, we here employ the \code{ProjectionOperator} provided by the TIGRE plugin, though the ASTRA plugin could equally have been used.

\begin{figure}[t]
\centering
\includegraphics[width=\linewidth]{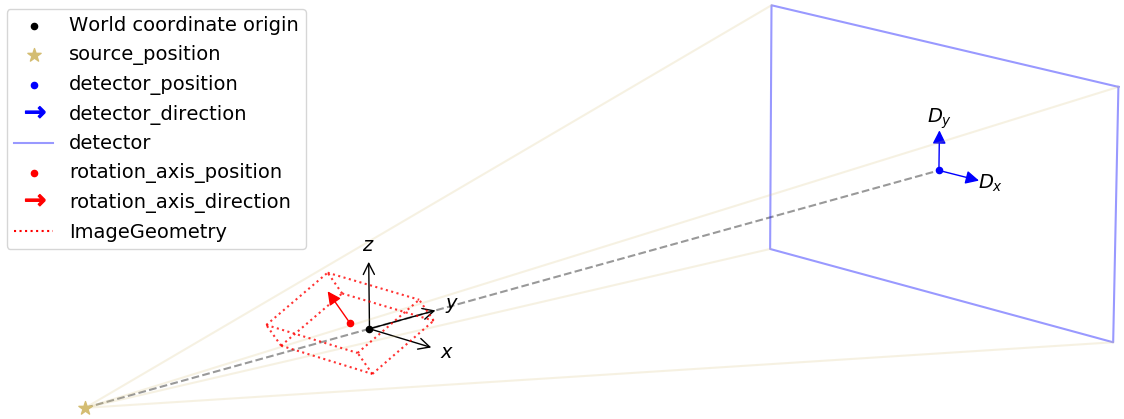}
\caption{CIL \code{AcquisitionGeometry} and \code{ImageGeometry} illustrated for the laminography cone-beam setup. Configurable parameters are shown in the legend. Parallel-beam geometry and 2D versions are also available. 
CIL can illustrate \code{ImageGeometry} and \code{AcquisitionGeometry} instances as in this figure using \code{show\_geometry}.
\label{fig:CILgeometrylami}}
\end{figure}

Laminography is an imaging technique designed for planar samples in which the rotation axis is tilted relative to the beam direction. Conventional imaging of planar samples often leads to severe limited-angle artifacts due to lack of transmission in-plane, while laminography can provide a more uniform exposure \cite{Xu:12}. In Transmission Electron Microscopy (TEM) the same technique is known as conical tilt.

An experimental laminography setup in the so-called rotary configuration was developed \cite{Fisher_2019} for Nikon micro-CT scanners in the Manchester X-ray Imaging Facility. Promising reconstructions of a planar LEGO-brick test phantom were obtained using the CGLS algorithm. Here we use CIL on the same data \cite{fisher_sarah_l_2019_2540509} to demonstrate how TV-regularisation and non-negativity constraints can reduce inherent laminographic reconstruction artifacts. 
CIL allows the specification of very flexible scan configurations. The cone-beam laminography setup of the LEGO data set provides an illustrative case for demonstrating CIL geometry, see \cref{fig:CILgeometrylami}. This particular geometry can be specified as follows, illustrating how different geometry components are used:
\begin{center}
\begin{tcolorbox}[
    enhanced,
    attach boxed title to top center={yshift=-2mm},
    colback=darkspringgreen!20,
    colframe=darkspringgreen,
    colbacktitle=darkspringgreen,
    title=Specify non-standard laminography acquisition geometry with full flexibility,
    text width = 15cm,
    fonttitle=\bfseries\color{white},
    boxed title style={size=small,colframe=darkspringgreen,sharp corners},
    sharp corners,
]
\begin{minted}{python}
ag = AcquisitionGeometry.create_Cone3D(                             \
        source_position=[0.0, source_pos_y, 0.0],                   \
        detector_position=[0.0, detector_pos_y, 0.0],               \
        rotation_axis_position=[object_offset_x, 0.0, 0.0],         \
        rotation_axis_direction=[0.0, -np.sin(tilt), np.cos(tilt)]) \
    .set_angles(angles=angles_list, angle_unit='degree')            \
    .set_panel(num_pixels=[num_pixels_x, num_pixels_y],             \
               pixel_size=pixel_size_xy, origin='top-left')   
\end{minted}
\end{tcolorbox}
\end{center}

\begin{figure}[tb]
\centering
\includegraphics[width=0.495\linewidth]{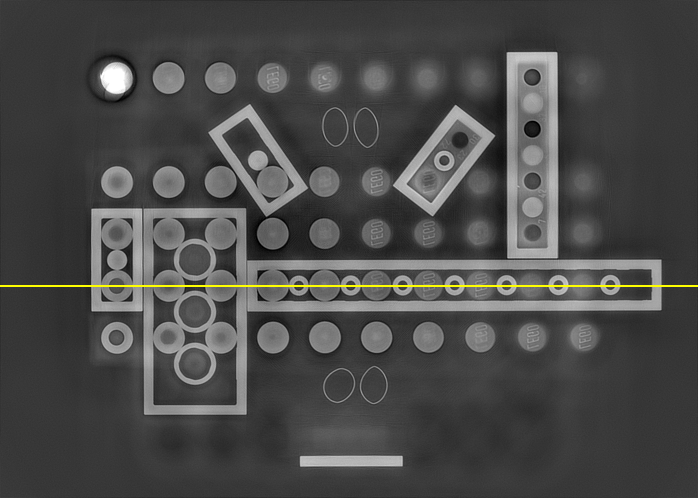}
\includegraphics[width=0.495\linewidth]{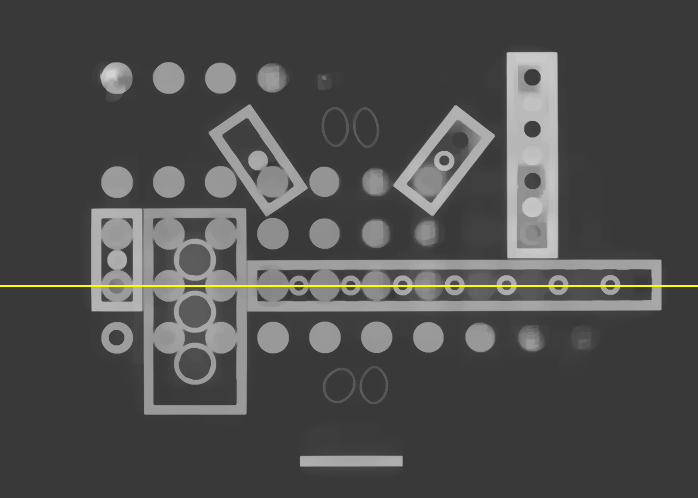}\\[0.5mm]
\includegraphics[width=0.495\linewidth]{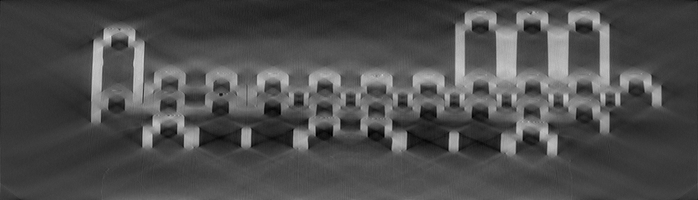}
\includegraphics[width=0.495\linewidth]{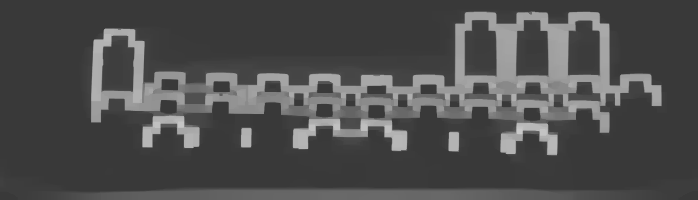} 
\caption{Slices through 3D reconstruction of laminography LEGO sample. Left, top/bottom: LS reconstruction using FISTA, horizontal/vertical slice at yellow line. Right: 
Same using TVNN, in which laminography artifact are suppressed while edges preserved.}
\label{fig:lami}
\end{figure}

The data consists of 2512 projections of 798-by-574 pixels sized 0.508 mm in a 360$^\circ$ cone-beam geometry.
We load the data with \code{NikonDataReader} and preprocess with a couple of \proc{} instances to prepare it for reconstruction. For reconstruction we use the GPU-accelerated cone-beam \code{ProjectionOperator} from \module{ccpi.plugin.tigre} and \code{FISTA} to solve \cref{eq:fistaproblem} for the unregularised least-squares problem (LS) and non-negativity constrained TV-regularised least-squares (TVNN). For TVNN we use \code{FBP_TV} from \module{cil.plugins.ccpi\_regularisation} which implements a GPU-accelerated version of $g_\text{TV}$, which is faster than, but otherwise equivalent to, using the native CIL \code{TotalVariation}. The full 3D volume is reconstructed for LS and TVNN, and \cref{fig:lami} shows a horizontal and vertical slice through both.

The LEGO bricks are clearly visualised in all reconstructions. The LS reconstruction has a haze in the horizontal slice (top left), which in the vertical slice (bottom left) is seen to amount to smooth directional streaks known to be inherent for laminography; in particular horizontal edges are heavily blurred. On the other hand, fine details in the horizontal plane are preserved, for example the text ``LEGO'' seen on several knobs to the right.

TVNN (right) reduces the haze and streaks substantially with the LEGO bricks displaying a uniform gray level and the horizontal edges in the vertical slice completely well-defined. However, some fine details are lost, including the ``LEGO'' text, which is a commonly observed drawback of TV-regularisation. Depending on the sample and application, this may or may not be an issue, and if necessary more sophisticated regularisers such as Total Generalised Variation (TGV) could be explored (a CIL example with TGV is given in the Part II article \cite{CIL2}). 

As shown, CIL can process very general scan configurations and allows easy experimentation with different reconstruction methods, including using third-party software through plugins.

\subsection{PET reconstruction in CIL using SIRF} \label{sec:exsirf}

SIRF (Synergistic Image Reconstruction Framework) \cite{SIRF2020} is an open-source platform for joint reconstruction of PET  and MRI  data developed by CCP-SyneRBI (formerly CCP-PETMR). CIL and SIRF have been developed with a large degree of interoperability, in particular data structures are aligned to enable CIL algorithms to work directly on SIRF data. As an example we demonstrate here reconstruction of the NEMA IQ Phantom \cite{NEMAstandard}, which is a standard phantom for testing scanner and reconstruction performance. It consists of a Perspex container with inserts of different-sized spheres, some filled with liquid with higher radioactivity concentration than the background, others with ``cold'' water (see \cite{NEMAstandard} for more details). This allows assessment of resolution and quantification. 

A 60-minute PET dataset \cite{thomas_benjamin_aneurin_2018_1304454} of the NEMA IQ phantom was acquired on a Siemens Biograph mMR PET/MR scanner at Institute of Nuclear Medicine, UCLH, London. Due to poor data statistics in PET a Poisson noise model is normally adopted, which leads to using the Kullback-Leibler (KL) divergence as data fidelity. We compare here reconstruction using the Ordered Subset Expectation Maximisation (OSEM) method \cite{OSEM} available in SIRF without using CIL, and TV-regularised KL divergence minimisation using CIL's \code{PDHG} algorithm with a \code{KullbackLeibler} data fidelity (KLTV). Instead of a CIL \op{} a SIRF \code{AcquisitionModel} represents the forward model, and has all necessary methods to allow its use in CIL algorithms. 

\begin{figure}[tb]
\includegraphics[width=0.305\linewidth,clip,trim=5.5cm 7.7cm 6.5cm 7cm]{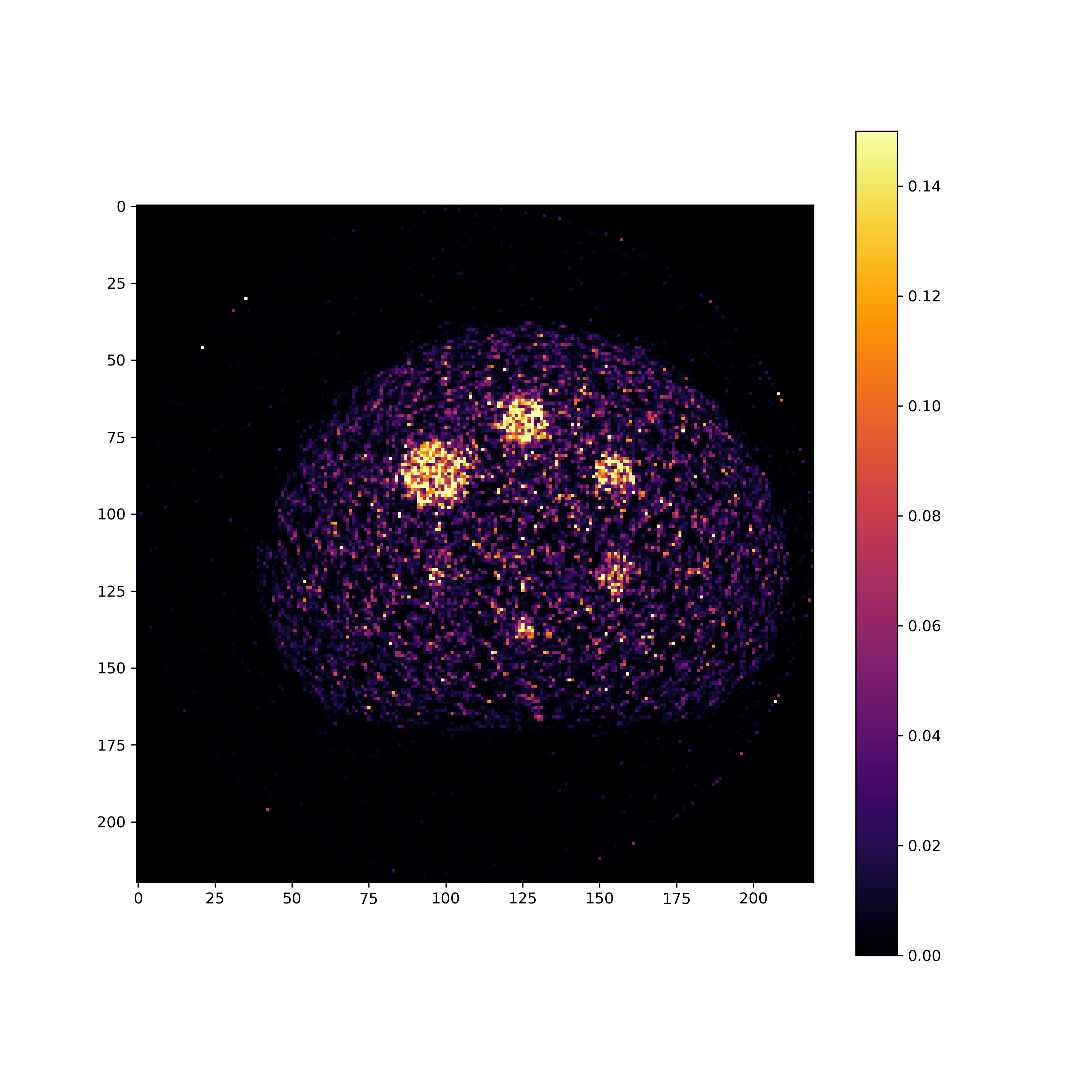}
\includegraphics[width=0.305\linewidth,clip,trim=5.5cm 7.7cm 6.5cm 7cm]{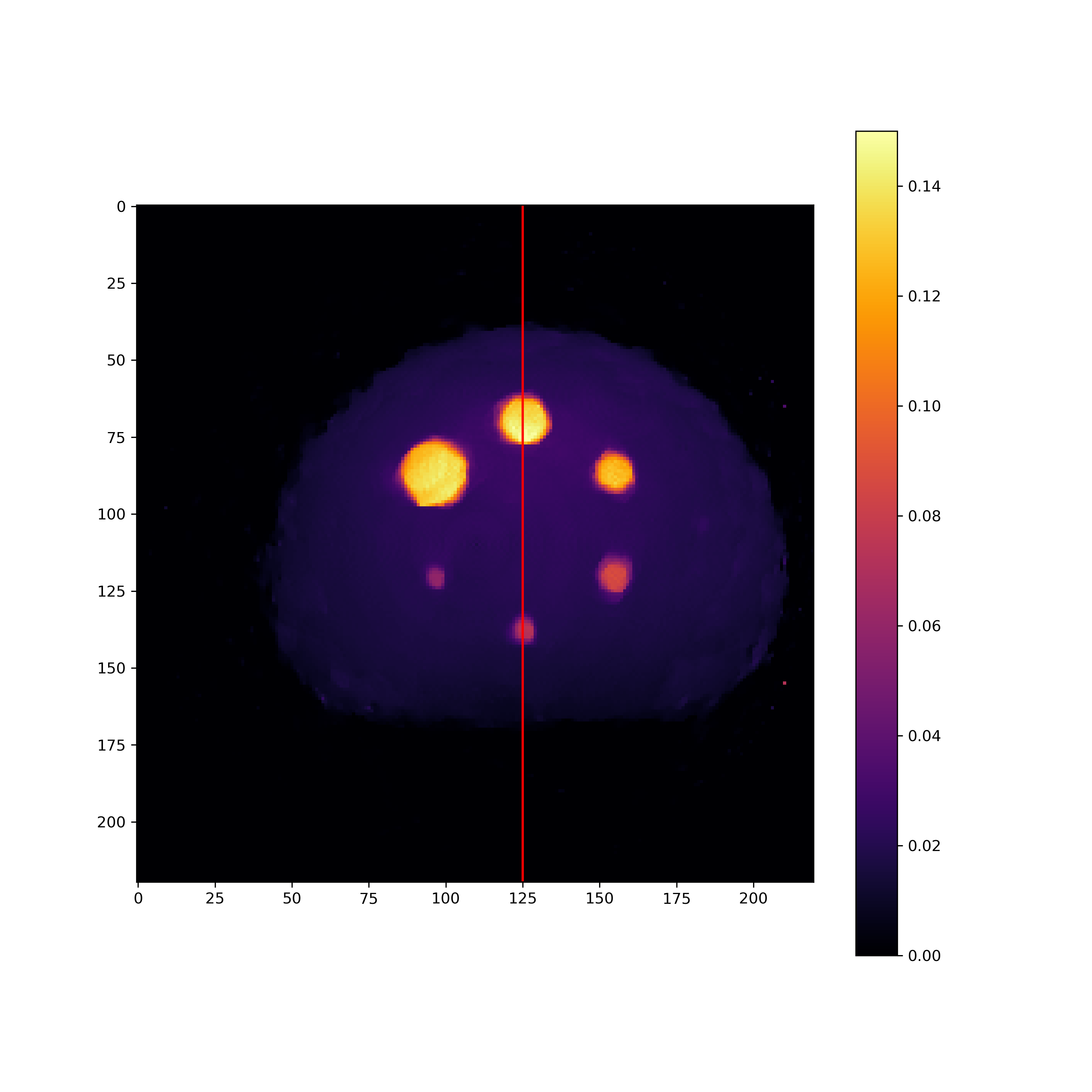}
\includegraphics[width=0.38\linewidth,clip,trim=0.5cm 0.7cm 1.5cm 1.2cm]{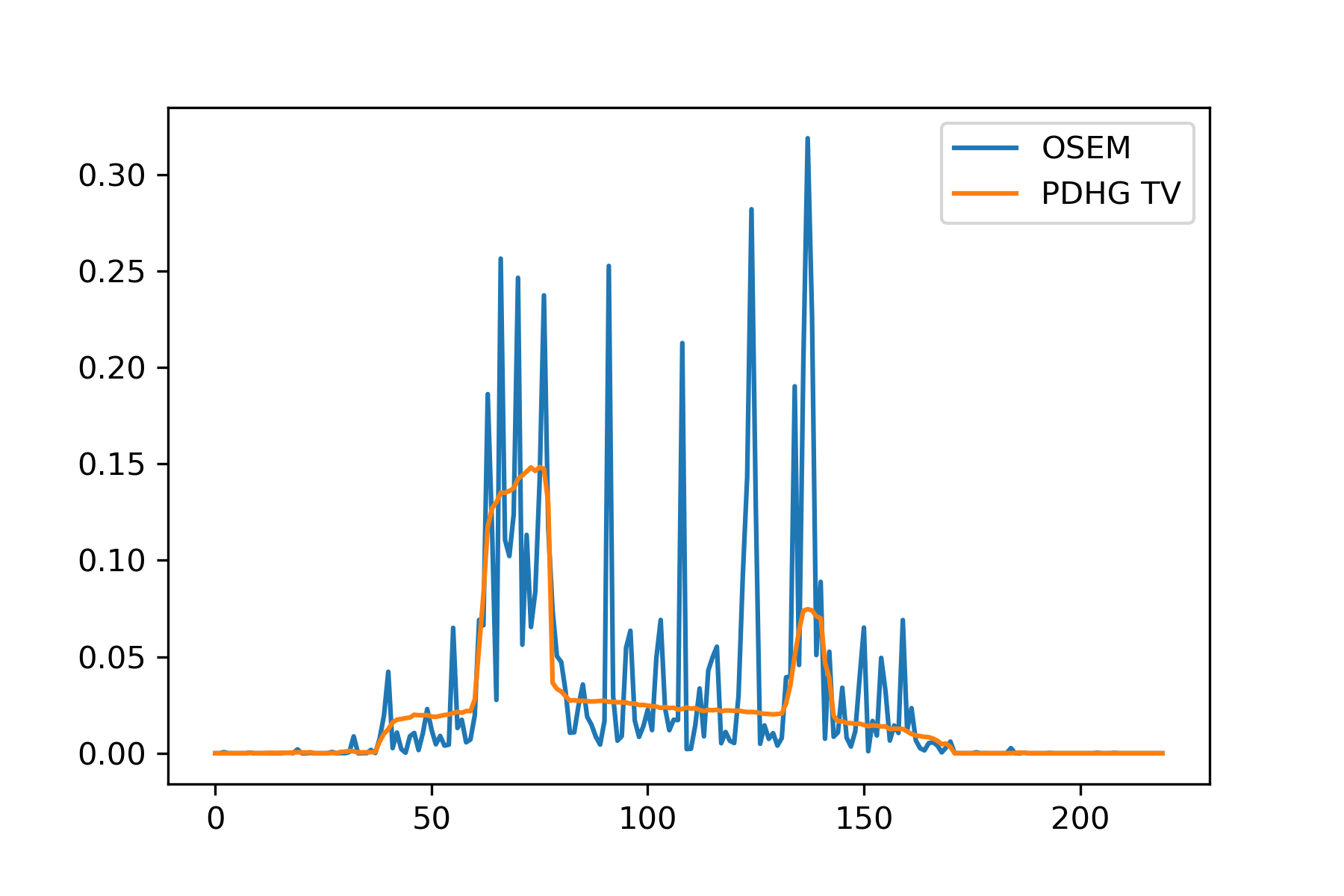}
\caption{3D PET reconstruction of NEMA IQ phantom data using CIL with SIRF data structures. Left: OSEM reconstruction (SIRF), horizontal slice. Centre: KLTV reconstruction (CIL PDHG). Colour range both [0,0.15]. Right: OSEM and KLTV profiles along red vertical line on centre plot. \label{fig:SIRFreduced}}
\end{figure}

\Cref{fig:SIRFreduced} shows horizontal slices through the $220\times 220\times 127$-voxel OSEM and KLTV reconstructions and vertical profile plots along the red line. In both cases the inserts are visible, but OSEM is highly affected by noise. KLTV reduces the noise dramatically, while preserving the insert and outer phantom edges. This may be beneficial in subsequent analysis, however a more detailed comparative study should take post-filtering into account. 

The purpose of this example was to give proof of 
principle of prototyping new reconstruction methods for PET with SIRF, using the generic algorithms of CIL, without needing to implement dedicated new algorithms in SIRF.
Another example with SIRF for PET/MR motion compensation employing CIL is given in \cite{SIRFspecialissue}.
 
\section{Summary and outlook} \label{sec:discussion}

We have described the CCPi Core Imaging Library (CIL), an open-source library, primarily written in Python, for processing tomographic data, with particular emphasis on enabling a variety of regularised reconstruction methods. The structure is highly modular to allow the user to easily prototype and solve new problem formulations that improve reconstructions in cases with incomplete or low-quality data. We have demonstrated the capability and flexibility of CIL across a number of test cases, including parallel-beam, cone-beam, non-standard (laminography) scan geometry, neutron tomography and PET using SIRF data structures in CIL. Further multi-channel cases including temporal/dynamic and spectral tomography are given in \cite{CIL2}.

CIL remains under active development with continual new functionality being added, steered by ongoing and future scientific projects. Current plans include:
\begin{itemize}
\item adding more algorithms, functions, and operators to support an even greater set of problems, for example allow convex constraints in smooth problems;
\item adding more pre-/postprocessing tools, for example to handle beam hardening;
\item adding templates with preselected functions, algorithms, etc. to simplify solving common problems such as TV regularisation;
\item further integrating with other third-party open-source tomography software through the plugin capability;
\item introducing support for nonlinear problems, such as polarimetric neutron spin tomography \cite{desai2020polarimetric} and electron strain tomography \cite{tovey2020scanning}; and
\item developing support for multi-modality problems.
\end{itemize}

CIL is developed as open-source on GitHub, and questions, feature request and bug reports submitted as issues are welcomed. Alternatively, the developer team can be reached directly at \url{CCPI-DEVEL@jiscmail.ac.uk}. 
CIL is currently distributed through the Anaconda platform; in the future additional modes of distribution such as Docker images may be provided. Installation instructions, documentation and training material is available from \url{https://www.ccpi.ac.uk/cil} as well as at \cite{CILReleases}, as are GitHub repositories with source code that may be cloned/forked and built manually. In this way users may modify and contribute back to CIL. 

Finally we emphasize that a multitude of optimization and regularization methods exist beyond those currently implemented in CIL and demonstrated in the present article. Recent overviews are given for example by \cite{Arridge2019,benning_burger_2018,chambolle_pock_2016,iterativereview} with new problems and methods constantly being devised. CIL offers a modular platform to easily implement and explore such methods numerically as well as apply them directly in large-scale imaging applications.

\subsection*{Data accessibility}
CIL version \cilversion{} as presented here is available through Anaconda; installation instructions are at \url{https://www.ccpi.ac.uk/cil}. 
In addition, CIL v21.0 and subsequent releases are archived at \cite{CILReleases}. 
Python scripts to reproduce all results are available from \cite{jakob_sauer_jorgensen_2021_4744394}. 
The steel-wire data set is provided as part of CIL; the original data is at \cite{DLSdata}. 
The neutron data set is available from \cite{jorgensenIMATzenodo}. 
The laminography data set is available from \cite{fisher_sarah_l_2019_2540509}.
The NEMA IQ PET data set is available from \cite{thomas_benjamin_aneurin_2018_1304454}.

\subsection*{Author contributions}
JJ designed and coordinated the study, carried out the steel-wire and neutron case studies, wrote the manuscript, and conceived of and developed the CIL software.
EA processed and analysed data for the neutron case study and developed the CIL software.
GB co-designed, acquired, processed and analysed the neutron data.
GF carried out the laminography case study and developed the CIL software.
EPap carried out the PET case study and developed the CIL software.
EPas conceived of and developed the CIL software and interoperability with SIRF.
KT contributed to the PET case study, interoperability with SIRF and  development of the CIL software.
RW assisted with case studies and contributed to the CIL software.
MT, WL and PW helped conceptualise and roll out the CIL software.
All authors critically revised the manuscript, gave final approval for publication and agree to be held accountable for the work performed therein.

\subsection*{Competing interests}
The authors declare that they have no competing interests.

\subsection*{Funding}
The work presented here was funded by the EPSRC grants ``A Reconstruction Toolkit for
Multichannel CT'' (EP/P02226X/1), ``CCPi: Collaborative Computational Project in Tomographic Imaging'' (EP/M022498/1 and EP/T026677/1), ``CCP PET-MR: Computational Collaborative Project in Synergistic PET-MR Reconstruction'' (EP/M022587/1) and ``CCP SyneRBI: Computational Collaborative Project in Synergistic Reconstruction for Biomedical Imaging'' (EP/T026693/1).
We acknowledge the EPSRC for funding the Henry Moseley X-ray Imaging Facility through grants (EP/F007906/1, EP/F001452/1, EP/I02249X/1, EP/M010619/1, and EP/F028431/1) which is part of the Henry Royce Institute for Advanced Materials funded by EP/R00661X/1. 
JSJ was partially supported by The Villum Foundation (grant no. 25893). 
EA was partially funded by the Federal Ministry of Education and Research (BMBF) and the Baden-W\"{u}rttemberg Ministry of Science as part of the Excellence Strategy of the German Federal and State Governments.
WRBL acknowledges support from
a Royal Society Wolfson Research Merit Award. 
PJW and RW acknowledge support from the European Research Council grant No. 695638 CORREL-CT. 
We thank Diamond Light Source for access to beamline I13-2 (MT9396) that contributed to the results presented here, and Alison Davenport and her team for the sample preparation and experimental method employed. 
We gratefully acknowledge beamtime RB1820541 at the IMAT Beamline of the ISIS Neutron and Muon Source, Harwell, UK. 

\subsection*{Acknowledgements}
We are grateful to input from Daniil Kazantsev for early-stage contributions to this work. We are grateful to Josef Lewis for building the neutron experiment aluminium sample holder and help with sample preparation at IMAT. We wish to express our gratitude to numerous people in the tomography community for valuable input that helped shape this work, including
Mark Basham,
Julia Behnsen,
Ander Biguri,
Richard Brown,
Sophia Coban,
Melisande Croft,
Claire Delplancke,
Matthias Ehrhardt,
Llion Evans,
Anna Fedrigo,
Sarah Fisher,
Parmesh Gajjar,
Joe Kelleher,
Winfried Kochelmann,
Thomas Kulhanek,
Alexander Liptak,
Tristan Lowe,
Srikanth Nagella,
Evgueni Ovtchinnikov,
S\o{}ren Schmidt,
Daniel Sykes,
Anton Tremsin,
Nicola Wadeson,
Ying Wang,
Jason Warnett, 
and 
Erica Yang. 
This work made use of computational support by CoSeC, the Computational Science Centre for Research Communities, through CCPi and CCP-SyneRBI.
	
	\bibliographystyle{vancouver}
\bibliography{CIL1paper.bib}

\end{document}